\documentclass[a4paper,11pt]{article}
\usepackage{latexsym}
\usepackage{amsmath}
\usepackage{amsthm}

\usepackage{graphicx}
\usepackage{txfonts}
\usepackage{bm}
\usepackage{color}
\usepackage{epic,eepic}

\usepackage{geometry}
\numberwithin{equation}{section}

\newtheorem{theorem}{Theorem}[section]
\newtheorem{proposition}[theorem]{Proposition}
\newtheorem{corollary}[theorem]{Corollary}
\newtheorem{lemma}[theorem]{Lemma}
\newtheorem{claim}[theorem]{Claim}
\newtheorem{remark}{Remark}[section]
\newtheorem{example}{Example}[section]

\newcommand{\OMIT}[1]{{\bf [OMIT:} #1 \ {\bf --- end OMIT] }}  
   \renewcommand{\OMIT}[1]{}            

\newcommand{\RR}{{\bf R}}
\newcommand{\ZZ}{{\bf Z}}

\newcommand{\finbox}{\hspace*{\fill}$\rule{0.17cm}{0.17cm}$}
\newcommand{\finboxHere}{\ $\rule{0.17cm}{0.17cm}$}

\newcommand{\odotZ}{\overset{....}} 

\newcommand{\llceil}{\bigg\lceil} 
\newcommand{\rrceil}{\bigg\rceil}

\newcommand{\llfloor}{\bigg\lfloor} 
\newcommand{\rrfloor}{\bigg\rfloor}

\newcommand{\Proof}{\noindent {\bf Proof.  }}

\usepackage{lineno}
\newcommand*\patchAmsMathEnvironmentForLineno[1]{
  \expandafter\let\csname old#1\expandafter\endcsname\csname #1\endcsname
  \expandafter\let\csname oldend#1\expandafter\endcsname\csname end#1\endcsname
  \renewenvironment{#1}
     {\linenomath\csname old#1\endcsname}
     {\csname oldend#1\endcsname\endlinenomath}}
\newcommand*\patchBothAmsMathEnvironmentsForLineno[1]{
  \patchAmsMathEnvironmentForLineno{#1}
  \patchAmsMathEnvironmentForLineno{#1*}}
\AtBeginDocument{
\patchBothAmsMathEnvironmentsForLineno{equation}
\patchBothAmsMathEnvironmentsForLineno{align}
\patchBothAmsMathEnvironmentsForLineno{flalign}
\patchBothAmsMathEnvironmentsForLineno{alignat}
\patchBothAmsMathEnvironmentsForLineno{gather}
\patchBothAmsMathEnvironmentsForLineno{multline}
}

\begin{document}

\title{
Decreasing Minimization on M-convex Sets:
\\
Background and Structures
}

\author{Andr\'as Frank\thanks{MTA-ELTE Egerv\'ary Research Group,
Department of Operations Research, E\"otv\"os University, P\'azm\'any
P. s. 1/c, Budapest, Hungary, H-1117. 
e-mail: {\tt frank@cs.elte.hu}. 
ORCID: 0000-0001-6161-4848.
The research was partially supported by the
National Research, Development and Innovation Fund of Hungary
(FK\_18) -- No. NKFI-128673.
}
\ \ and \ 
{Kazuo Murota\thanks{Department of Economics and Business Administration,
Tokyo Metropolitan University, Tokyo 192-0397, Japan, 
e-mail: {\tt murota@tmu.ac.jp}. 
Currently at
The Institute of Statistical Mathematics,
Tokyo 190-8562, Japan.
ORCID: 0000-0003-1518-9152.
The research was supported by CREST, JST, Grant Number JPMJCR14D2, Japan, 
and JSPS KAKENHI Grant Numbers JP26280004, JP20K11697. 
 }}}

\date{July 2020 / July 2021} 


\maketitle

\begin{abstract} 
The present work is the first member of a pair of papers concerning
decreasingly-minimal (dec-min) elements of a set of integral vectors,
where a vector is dec-min if its largest component is as small as
possible, within this, the next largest component is as small as
possible, and so on.  
This discrete notion, along with its fractional counterpart, showed
up earlier in the literature under various names.

The domain we consider is an M-convex set, that is, 
the set of integral elements of an integral base-polyhedron.  
A fundamental difference between the fractional and the discrete case is that a
base-polyhedron has always a unique dec-min element, while 
the set of dec-min 
elements of an M-convex set 
admits a rich structure, described
here with the help of a \lq canonical chain\rq.  As a consequence, we
prove that this set arises from a matroid by translating the
characteristic vectors of its bases with an integral vector.

By relying on these characterizations, we prove that an element is
dec-min if and only if the square-sum of its components is minimum, a
property resulting in a new type of min-max theorems.  
The characterizations also give rise, as shown in the companion paper,
to a strongly polynomial algorithm, and to several applications in the areas
of resource allocation, network flow, matroid, and graph orientation problems, 
which actually provided a major motivation to the present investigations.  
In particular, we prove a conjecture on graph orientation. 
\end{abstract}

{\bf Keywords}:  Submodular optimization, Matroid, Base-polyhedron,
M-convex set, Lexicographic minimization.

{\bf Mathematics Subject Classification (2010)}: 90C27, 05C, 68R10


{\bf Running head}:
Decreasing Minimization: Background and Structures

\newpage
\tableofcontents

\newpage

\section{Introduction}
\label{SCintro}

We investigate a problem which we call ``discrete decreasing minimization.''  
An element of a set of vectors is called
decreasingly minimal (dec-min) if its largest component is as small as possible, 
within this, its second largest component is as small as possible, and so on.  
The term discrete decreasing minimization refers
to the problem of finding a dec-min element 
(or even a cheapest dec-min element with respect to a given weighting)
of a set of integral vectors.
In the present work, this set is an M-convex set, which is nothing but
the set of integral elements of an integral base-polyhedron.  
Note that one may consider the analogous term ``increasing maximization'' (inc-max), as well.
This dichotomy is the reason why we avoid the usage of term 
``lexicographic optimal'' used in the literature.

The goal of this paper is to develop structural characterizations of the
set of dec-min elements of an M-convex set.  
These form the bases, in \cite{FM19partB}, 
for developing a strongly polynomial algorithm, 
as well as for exploring and exhibiting various applications.  
Actually, earlier special cases played a major motivating role for our investigations,
which are described in Section~\ref{SCbackprob}.  
The main results will be described in Section~\ref{SCmaingoal}.
The research was strongly motivated by the theory of 
Discrete Convex Analysis (DCA),
but the paper is self-contained and does not rely on any prerequisite from DCA.

\subsection{Background problems}
\label{SCbackprob}

There are several different sources 
underlying the study of discrete decreasing minimization.

\subsubsection{Orientations of graphs}
\label{SCorigraph}

\begin{figure}[b]
\centering
\includegraphics[width=0.8\textwidth,clip]{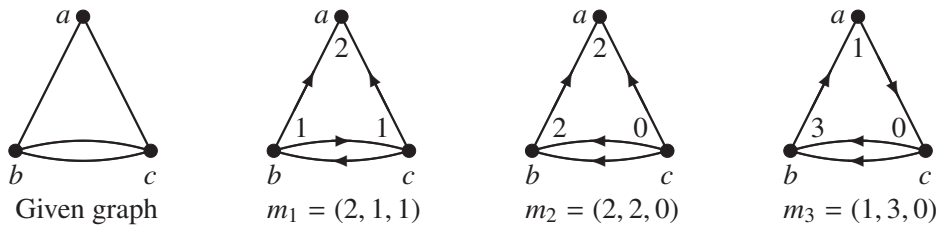}
\caption{Orientations of a graph} 
\label{FGoritri}
\end{figure}

Let $G=(V,E)$ be an undirected graph.  
Orienting an edge $e=uv$ means the operation that replaces $e$ 
by one of the two oppositely directed edges 
(sometimes called arcs) $uv$ or $vu$.  
A directed graph arising from $G$ by orienting all of its edges 
is called an orientation of $G$ (see Figure \ref{FGoritri}).  
A graph orientation problem consists of finding an orientation
of $G$ meeting some specified properties such as in-degree constraints
(lower and upper bounds) and/or various connectivity prescriptions.
One goal is to characterize undirected graphs for which the requested
orientation exists, and 
another related goal is to design an algorithm for finding the orientation.

The literature is quite rich in orientation results,
for a relatively wide overview, see the book \cite{Frank-book}.  
While upper and lower bounds 
are often imposed on the in-degree of each node in an orientation of $G$,
there are other type of requirements 
for an orientation concerning the global distribution of the in-degrees of nodes.  
That is, the goal is to find orientations 
(with possible connectivity expectations) 
whose in-degree vector (on the node-set) is felt intuitively evenly distributed:  
\lq fair\rq, \lq equitable\rq, \lq egalitarian\rq,  \lq uniform\rq.  
For example, how can one determine 
the minimum value $\beta_{1}$ of the largest
in-degree of a ($k$-edge-connected) orientation?  
Even further, after determining $\beta_{1}$, it might
be interesting to minimize the number of nodes 
with in-degree $\beta_{1}$ among orientations of $G$ 
with largest in-degree $\beta_{1}$.  
Or, a more global equitability feeling is captured if we minimize the sum
of squares of the in-degrees.  
For example, the in-degree vector $m_{3}=(1,3,0)$
in Figure \ref{FGoritri} (with square-sum 10) 
is felt less \lq fair\rq\ 
than $m_{2}=(2,2,0)$ (with square-sum 8) and 
$m_{1}=(2,1,1)$ (with square-sum 6).

A formally different definition was recently suggested and
investigated by Borradaile et al.~\cite{BIMOWZ} who called an
orientation of $G=(V,E)$ {\bf egalitarian} if the highest in-degree of
the nodes is as small as possible, and within this, the second highest
(but not necessarily distinct) in-degree is as small as possible, and
within this, the third highest in-degree is as small as possible, and
so on.  In other words, if we rearrange the in-degrees of nodes in a
decreasing order, then the sequence is lexicographically minimal.  
In order to emphasize that the in-degrees are considered in a decreasing
order, we prefer the use of the more expressive term 
{\bf decreasingly minimal} ({\bf dec-min}, for short) 
for such an orientation, rather than egalitarian.

This change of terminology is reasonable since one may also consider
the mirror problem of finding an {\bf increasingly maximal}
 (or {\bf inc-max}, for short) orientation that is an orientation of $G$ 
in which \ the smallest in-degree is as large as possible, within this,
the second smallest in-degree is as large as possible, and so on.
Intuitively, such an orientation may equally 
be felt \lq egalitarian\rq \ in the informal meaning of the word.

Borradaile et al.~\cite{BIMOWZ}, however, proved that 
an orientation of a graph is
decreasingly minimal (egalitarian in their original term) 
if and only if 
there is no \lq small\rq \ improvement, where a small improvement
means the reorientation of a dipath from some node $s$ to another node
$t$ with in-degrees $\varrho(t)\geq \varrho(s)+2$.
This theorem immediately implies that an orientation is decreasingly minimal 
if and only if 
it is increasingly maximal, and therefore we could retain the
original terminology ``egalitarian orientation'' used in \cite{BIMOWZ}.

However, when orientations are considered with specific requirements
such as strong (or, more generally, $k$-edge-) connectivity and/or
in-degree bounds on the nodes, the possible equivalence of
decreasingly minimal and increasingly maximal orientations had not yet been investigated.  
Actually, Borradaile et al.~\cite{BIMOWZ}
conjectured that a strong orientation of a graph 
is decreasingly minimal (among strong orientations) 
if and only if 
there is no small improvement preserving
strong connectivity.
This, if true, would imply immediately that
decreasing minimality and increasing maximality do coincide for strong
orientations, as well.
In \cite{FM19partB} we shall prove this conjecture even in its extended form 
concerning $k$-edge-connected and in-degree constrained orientations.

\subsubsection{A resource allocation problem}
\label{SCresource}

Another source of our investigations is due to Harvey et al.~\cite{HLLT} 
who solved the problem of minimizing 
$\sum [d_{F}(s) (d_{F}(s)+1):  s\in S]$ 
over the semi-matchings $F$ of a simple bipartite graph $G=(S,T;E)$.  
Here a semi-matching is a subset $F$ of edges 
for which $d_{F}(t)=1$ holds for every node $t\in T$.  
Harada et al.~\cite{HOSY07} solved the minimum 
edge-cost version of this problem.
The framework of Harvey et al. was extended by Bokal et al.~\cite{BBJ} 
to quasi-matchings, and, even further, to degree-bounded quasi-matchings
by Katreni{\v c} and Semani{\v s}in \cite{Katrenic13}.
It turns out that these problems are strongly related to minimization 
of a separable convex function over (integral elements of) 
a base-polyhedron which has been investigated in the literature 
under the name of ``resource allocation problems under submodular constraints'' 
(\cite{Fed-Gro86}, \cite{Hochbaum-Hong}, \cite{Hoc07},
\cite{KSI}, \cite{Ibaraki-Katoh.88}, \cite{Katoh-Ibaraki.98}).  
Ghodsi et al.~\cite{GZSS} considered the problem of finding a semi-matching
$F$ of $G=(S,T;E)$ whose degree-vector restricted to $S$ is
increasingly maximal.  
This problem ``constrained max-min fairness'' 
originated from modelling a fair
sharing problem for datacenter jobs.  
Here $T$ corresponds to the set of available computers while $S$ to the set of users.  
An edge $st$ exists
if user $s$ can run her program on computer $t$.
Ghodsi et al. also consider the fractional version, which finds a real vector
$x:E\rightarrow {\bf R}_{+}$ so that $d_{x}(t)=1$ for every $t\in T$ and
the vector $(d_{x}(s):  s\in S)$ is increasingly maximal. 
 (Here $d_{x}(v):= \sum [ x(uv):  uv\in E]$).  
When $x$ is requested to be $\{ 0,1 \}$-valued, we are back at the subgraph version.
It should be emphasized that, 
unlike the well-known situation with ordinary bipartite matchings,
the optima for the subgraph version and for the fractional version 
may be different.

\subsubsection{Network flows}
\label{SCnetwork}

There is a much earlier, strongly related problem concerning network flows, 
due to Megiddo \cite{Megiddo74}, \cite{Megiddo77}.  
We are given a digraph $D=(V,A)$ with a source-set $S\subset V$ 
and a single sink-node $t\in V-S$. 
Let $g:A\rightarrow {\bf R}_{+}$ be a capacity function.  
A flow means a function 
$x:A\rightarrow {\bf R}_{+}$ 
for which the net out-flow
$\delta_{x}(v)-\varrho_{x}(v)=0$ if $v\in V-(S \cup \{ t \})$ and 
$\delta_{x}(v)-\varrho_{x}(v)\geq 0$ if $v\in S$.
 (Here $\varrho_{x}(v):= \sum [x(uv):uv\in A]$ and 
$\delta_{x}(v):= \sum [x(vu):  vu\in A]$.) 
A flow $x$ is feasible if $x\leq g$.  
The {\bf flow  amount} of $x$ is defined by $\varrho_{x}(t)-\delta_{x}(t)$.
Megiddo solved the problem of finding a feasible flow of maximum flow amount which is,
in his term, ``source-optimal'' at $S$.  
Source-optimality is the same as requiring that the net out-flow vector 
on $S$ is increasingly maximal.  
It must be emphasized that the flow in Megiddo's problem 
is not requested to be integer-valued.

The integrality property is a fundamental feature of ordinary network flows.  
It states that in case of an integer-valued capacity function
$g$ there always exists a maximum flow which is integer-valued.  
In this light, it is quite surprising that the integer-valued
 (or discrete) version of Megiddo's inc-max problem
(source-optimal in his term), 
where the capacity function $g$ is integer-valued 
and the max flow is required to be integer-valued, has not been investigated
in the literature.
We consider the present work the first such attempt.

\subsubsection{Matroid bases}

The fourth source of discrete decreasing minimization problems is due to
Levin and Onn \cite{Levin-Onn} who used the term ``shifted optimization.''  
They considered the following matroid optimization problem.  
For a specified integer $k$, find $k$ bases 
$Z_{1},Z_{2},\dots ,Z_{k}$ of a matroid $M$ on $S$ in such a way that the vector 
$\sum_{i}\chi_{Z_{i}}$ be, in our term, decreasingly minimal, where $\chi_{Z}$
is the incidence (or characteristic) vector of a subset $Z$.  
They apply the following natural approach to reduce the problem to classic
results of matroid theory.  
First replace each element $s$ of $S$ by
$k$ copies to be parallel in the resulting matroid $M'$ on the new ground-set 
$S'=S_{1}\cup S_{2}\cup \cdots \cup S_{k}$ 
where $S_{1},\dots ,S_{k}$
are the $k$ copies of $S$.  
Assign then a \lq rapidly increasing\rq \  cost function to the copies.  
(The paper \cite{Levin-Onn} explicitly
describes what rapidly increasing means).  
Then a minimum cost basis of the matroid $M_{0}$ obtained 
by multiplying $M'$ $k$-times will be a solution to the problem.  
(By definition, a basis of $M_{0}$
 is the union of $k$ disjoint bases of $M'$).

\subsection{Main goals}
\label{SCmaingoal}

Each of the problems in Section~\ref{SCbackprob} 
may be viewed as a special case of a single discrete optimization problem:  
\textit{Characterize decreasingly minimal elements of an M-convex set}
\cite{Murota98a,Murota03} 
(or, in other words, dec-min integral elements of an integral base-polyhedron).  
By one of its equivalent definitions, an M-convex set is nothing but 
the set of integral elements of an integral base-polyhedron.

We characterize dec-min elements of an M-convex set as those admitting
no local improvement, and prove that the set of dec-min elements is
itself an M-convex set arising by translating a matroid
base-polyhedron with an integral vector.  
This result implies that decreasing minimality and 
increasing maximality coincide for M-convex sets.  
We shall also show that an element of an M-convex set 
is dec-min precisely if it is a square-sum minimizer.  
Using the characterization of dec-min elements, we shall derive a novel min-max
theorem for the minimum square-sum of elements of an integral member
of a base-polyhedron.

The structural description of the set of dec-min elements of an M-convex set 
in terms of a matroid makes it possible 
to solve the algorithmic problem of finding a minimum cost dec-min element.  
(In the continuous counterpart of decreasing minimization, 
this problem simply does not exist due to the uniqueness
of the fractional dec-min element of a base-polyhedron.)  
In the companion paper \cite{FM19partB}, we shall 
describe a polynomial algorithm
for finding a minimum cost (in-degree constrained) dec-min orientation.  
Furthermore, we shall outline an algorithm to solve the
minimum cost version of the resource allocation problem of 
Harvey et al.~\cite{HLLT} mentioned in Section~\ref{SCresource}.
Moreover, as an essential extension of the algorithm of Harada et al.~\cite{HOSY07}, 
we describe a strongly polynomial algorithm to
solve a minimum cost version of 
the decreasingly minimal degree-bounded subgraph problem 
in a bipartite graph $G=(S,T,E)$.
One may consider two versions here.  
In the simpler one, we have a cost-function on the node-set of $G$, 
that is, on the ground-set of the corresponding M-convex set.  
Due to the matroidal description of the set of dec-min elements of an M-convex set, 
this min-cost version becomes rather easy since the matroid greedy algorithm can be applied.
Significantly more complicated is the other min-cost version,
where a cost-function is given on the set of edges of $G$.
The latter version is also solved in \cite{FM19partB}.

The topic of our investigations may be interpreted 
as a discrete counterpart of the work by Fujishige \cite{Fujishige80} from 1980 
on the lexicographically optimal base of a base-polyhedron $B$, where
lexicographical optimality is essentially the same as decreasing minimality.
He proved that there is a unique lexicographically optimal member $x_{0}$ of $B$, 
and $x_{0}$ is the unique minimum norm 
(that is, the minimum square-sum) element of $B$.  
This uniqueness result reflects a characteristic difference 
between the behaviour of the fractional and the discrete versions 
of decreasing minimization since
in the latter case the set of dec-min elements (of an M-convex set) is
typically not a singleton, and it actually has, as indicated above, 
a matroidal structure.  
While the present paper focuses on the unweighted case, 
the lexicographically optimal base of a base-polyhedron 
is defined and analyzed with respect to a weight vector in \cite{Fujishige80}.

Fujishige also introduced the concept of principal partitions
concerning the dual structure of the minimum norm point of a base-polyhedron.  
Actually, he introduced a special chain of the subsets of ground-set $S$ 
and his principal partition arises by taking the difference sets of this chain.  
We will prove that there is an analogous concept in the discrete case, as well.  
As an extension of the above-mentioned elegant result of Borradaile et al.~\cite{BMW}
concerning graphs, we show that there is a `canonical chain'
describing the structure of dec-min elements of an M-convex set.  
The relation between our canonical partition and
Fujishige's principal partition is clarified in \cite{FM19partII},
showing that the canonical partition is an intrinsic structure of an M-convex set 
consistent with the principal partition of a base-polyhedron.

The paper of Fujishige is one of the early representatives of the
rich literature of related topics.  
The section of ``Survey of early papers'' in \cite{FM19partII} 
provides a relatively complete overview of these results, 
along with an outline of their relationship.

The present paper is organized, as follows.  
After formally introducing the basic notions, terminology, and notation 
in Section~\ref{egal2base}, we prove, in Section~\ref{SCchardecmin},
two characterizations of an element $m$ of an M-convex sets to be
decreasingly minimal.
The first one is a co-{\bf NP}-characterization
consisting of an easily checkable certificate for $m$ not to be dec-min, 
while the second one is an {\bf NP}-characterization
consisting of an easily checkable certificate for $m$ to be dec-min.  
The first characterization implies immediately that 
$m$ is dec-min precisely if it is inc-max.  
This is a property that fails to hold for the intersection of two M-convex sets.

In Sections \ref{peak} and \ref{canonical}, 
we show how the set of all dec-min elements of an M-convex set 
can be obtained from a matroid, and hence it is also an M-convex set.  
The main device is the dual concept of canonical chains and partitions.
In Section~\ref{SCsqdiffsum}, we prove that 
an element $m$ of an M-convex set is dec-min if and only if its 
$\ell_{2}$-norm $\| \cdot \|_{2}$
(or equivalently, the square sum of its components) is minimum.  
On one hand, this result seems surprising 
in the light of the fact that 
there is an example for two elements $m_{1}$ and $m_{2}$ of an M-convex set 
for which $m_{1}$ is decreasingly smaller than $m_{2}$, but 
$\| m_{1} \|_{2} > \| m_{2} \|_{2}$.
On the other hand, we shall use the coincidence of
dec-minimality and square-sum minimality to prove a min-max theorem
for the minimum square-sum of the components 
over the elements of an M-convex set.  
As a special case, this provides a min-max formula for
the the minimum of 
$\sum \{\varrho_{D}(v)\sp{2} :  v\in V \}$ 
over all orientations $D$ of an undirected graph on node-set $V$, 
where $\varrho_{D}(v)$ denotes the in-degree of $v$ in $D$.
To our best knowledge, no min-max formulas of similar type appeared 
earlier in the literature.

In Section~\ref{SCcontRdiscZ}, 
we discuss the relationship between the continuous and discrete versions
of decreasing minimization on a base-polyhedron.
Both dec-min elements and partitions of the ground-set
are compared between the continuous and discrete versions.
Finally in Section~\ref{SCconcl}, 
we give a perspective of our series of research on discrete decreasing minimization.

\subsection{Notation}

Throughout the paper, $S$ denotes a finite non-empty ground-set.  
For elements $s,t\in S$, we say that $X\subset S$ is an {\bf $s\overline{t}$-set}
if $s\in X\subseteq S-t$.  
For a vector $m\in {\bf R}\sp{S}$ (or function $m:S\rightarrow {\bf R}$), 
the restriction of $m$ to $X\subseteq S$ is denoted by $m\vert X$.  
We also use the notation
$\widetilde m(X)=\sum [m(s):  s\in X]$.  
With a small abuse of notation, 
we do not distinguish between a one-element set $\{s\}$
called a singleton and its only element $s$.  
When we work with a chain 
$\cal C$ of non-empty sets $C_{1} \subset C_{2} \subset \cdots \subset C_{q}$, 
we sometimes use $C_{0}$ to denote the empty set without assuming
that $C_{0}$ is a member of ${\cal C}$.
The characteristic (or incidence) vector of a subset 
$Z$ is denoted by $\chi_{Z}$,  that is,  
$\chi_{Z}(s)=1$ if $s\in Z$ and $\chi_{Z}(s)=0$ otherwise.
For a polyhedron $B$, $\odotZ{B}$ 
 (pronounce:  dotted $B$) denotes the set of integral members 
(elements, vectors, points) of $B$, 
that is,
\begin{equation}    \label{odotBdef} 
 \odotZ{B} := B\cap \ZZ\sp{S}.
\end{equation}

For a set-function $h$, we allow it to have value $+\infty $ or $-\infty$,
while $h(\emptyset )=0$ is assumed throughout.
Where $h(S)$ is finite, the {\bf complementary function} $\overline{h}$ is
defined by $\overline{h}(X) = h(S) - h(S-X)$.  
For functions $f:S\rightarrow {\bf Z}\cup \{-\infty \}$ and
$g:S\rightarrow {\bf Z}\cup \{+\infty \}$ 
with $f\leq g$,  the polyhedron 
$T(f,g)=\{x\in {\bf R}\sp{S}:  f\leq x\leq g\}$
is called a {\bf box}.  
If $g(s)\leq f(s)+1$ holds for every $s\in S$, we speak of a 
{\bf small box}.  
For example, the $(0,1)$-box is small, and so is any set
consisting of a single integral vector.


\section{Base-polyhedra and M-convex sets} 
\label{egal2base}

Let $S$ be a finite non-empty ground-set.
Let $b$ be a set-function for which $b(X)=+\infty $ is allowed but
$b(X)=-\infty $ is not.  The submodular inequality for subsets
$X,Y\subseteq S$ is defined by 
\[
b(X) + b(Y) \geq b(X\cap Y) + b(X\cup Y).
\] 
We say that $b$ is {\bf submodular} if the submodular inequality holds
for every pair of subsets $X, Y\subseteq S$ with finite $b$-values.
A set-function $p$ is {\bf supermodular} if $-p$ is submodular.

For a submodular integer-valued set-function $b$ on $S$ 
for which $b(\emptyset )=0$ and $b(S)$ is finite,
the {\bf base-polyhedron} $B$ 
in ${\bf R}\sp{S}$ is defined by 
\begin{equation}    \label{Bsubfn}
B =B(b)
=\{x\in {\bf R}\sp{S}:  \widetilde x(S)=b(S), \ \widetilde x(Z)\leq b(Z)
  \ \mbox{ for every } \  Z\subset S\},
\end{equation}  
which is possibly unbounded.

A special base-polyhedron is the one of matroids.  
Given a matroid $M$, Edmonds proved that the polytope 
(that is, the convex hull) of the incidence 
(or characteristic) vectors of the bases of $M$ 
is the base-polyhedron $B(r)$ defined by the rank function $r$ of $M$, 
that is, 
$B(r)=\{x\in {\bf R}\sp{S}:  \ \widetilde x(S)=r(S)$ 
and
$\widetilde x(Z)\leq r(Z)$ for every subset $Z\subset S\}$. 
It can be proved that a kind of converse also holds, namely, every
(integral) base-polyhedron in the unit $(0,1)$-cube is a matroid base-polyhedron.
We call the translation of a matroid base-polyhedron 
a {\bf translated matroid base-polyhedron}.  
It follows that the intersection of a
base-polyhedron with a 
small box is a translated matroid base-polyhedron.

A base-polyhedron $B(b)$ is never empty, 
and $B(b)$ is known to be an integral polyhedron.  
(A rational polyhedron is {\bf integral} 
if each of its faces contains an integral element.  
In particular, a pointed rational polyhedron is integral 
if all of its vertices are integral.)
By convention, the empty set is also considered a base-polyhedron.
Note that a real-valued submodular function $b$ also defines a
base-polyhedron $B(b)$ but in the present work we are interested only
in integer-valued submodular functions and integral base-polyhedra.

We call the set $\odotZ{B}$ 
of integral elements of an integral base-polyhedron $B$ an {\bf M-convex set}.  
Originally, this basic notion of DCA 
introduced by Murota \cite{Murota98a} (see, also the book \cite{Murota03}), 
was defined as a set of integral points in 
${\bf R}\sp{S}$ satisfying certain exchange axioms, and it is known 
that the two properties are equivalent (\cite[Theorem 4.15]{Murota03}).  
While an integral base-polyhedron $B$ defines
an M-convex set as $\odotZ{B} = B \cap \ZZ\sp{S}$,
an M-convex set induces an integral base-polyhedron as its convex hull.  
This implies, in particular, that two distinct 
integral base-polyhedra $B_{1}$ and $B_{2}$
define distinct M-convex sets
$\odotZ{B}_{1}$ and $\odotZ{B}_{2}$.   
The set of integral elements of a translated matroid
base-polyhedron will be called a {\bf matroidal M-convex set}.

A non-empty base-polyhedron $B$ can also be defined by a supermodular
function $p$ for which $p(\emptyset )=0$ and $p(S)$ is finite as follows:  
\begin{equation}    \label{Bsupfn}
B=B'(p)=\{x \in {\bf R}\sp{S}:  
     \widetilde x(S)=p(S),  
     \widetilde x(Z)\geq p(Z) \ \mbox{ for every } \  Z\subset S\}.
\end{equation}  
It is known that $B$ uniquely determines both $p$ and $b$, namely,
$b(Z) = \max\{\widetilde x(Z):  x\in B\}$ 
and 
$p(Z) = \min\{\widetilde x(Z):  x\in B\}$.  
The functions $p$ and $b$ are complementary functions, 
that is, $b(X)= p(S) - p(S-X)$ or $p(X)= b(S) - b(S-X)$
(where $b(S)=p(S)$).

For a set $Z\subset S$, $p\vert Z$ 
denotes the restriction of $p$ to $Z$,
while $p'=p/Z$ is the set-function on $S-Z$ obtained from $p$ by contracting $Z$, which is
defined for $X\subseteq S-Z$ by $p'(X)= p(X\cup Z)-p(Z)$.  
,Note that $p/Z$ and $\overline{p}\vert (S-Z)$ are complementary set-functions.  
It is also known for disjoint subsets $Z_{1}$ and $Z_{2}$ of $S$ that
\begin{equation}
 (p/Z_{1})/Z_{2} = p/(Z_{1}\cup Z_{2}), 
\label{(Z1Z2)} 
\end{equation}
When $p(Z)$ is finite, the base-polyhedron $B'(p\vert Z)$
is called the {\bf restriction} of $B'(p)$ to $Z$.

Let $\{S_{1},\dots ,S_{q}\}$ be a partition of $S$ and let $p_{i}$ be a
supermodular function on $S_{i}$.  
Let $p$ denote the supermodular function on $S$ defined by 
$p(X):  = \sum [p_{i}(S_{i}\cap X):  i=1,\dots ,q]$ for $X\subseteq S$.  
The base-polyhedron $B'(p)$ is called the
{\bf direct sum} of the $q$ base-polyhedra $B'(p_{i})$.  
Obviously, a vector $x\in {\bf R}\sp{S}$ is in $B'(p)$ 
if and only if 
each $x_{i}$ is in $B'(p_{i})$ \ $(i=1,\dots ,q)$,
where $x_{i}$ denotes the restriction
$x\vert S_{i}$ of $x$ to $S_{i}$.

It is known that a face $F$ of a non-empty base-polyhedron $B$ is also
a base-polyhedron.  
The (special) face of $B'(p)$ defined by the single equality 
$\widetilde x(Z)=p(Z)$ is the direct sum of the base polyhedra 
$B'(p\vert Z)$ and $B'(p/Z)$.  
More generally, any face $F$ of $B$ 
can be described with the help of a chain 
$(\emptyset \subset ) \ C_{1} \subset C_{2} \subset \cdots \subset C_{\ell} = S$ 
of subsets by 
$F:= \{z:  z\in B, \ p(C_{i})=\widetilde z(C_{i})$ 
for $i=1,\dots ,\ell\}$.
(In particular, when $\ell=1$, the face $F$ is $B$ itself.)  
Let
$S_{1}:=C_{1}$ and $S_{i}:=C_{i}-C_{i-1}$ 
for $i=2,\dots ,\ell$.  
Then $F$ is the direct sum of the base-polyhedra $B'(p_{i})$,
where $p_{i}$ is a supermodular function on $S_{i}$ defined by 
$p_{i}(X):= p(X\cup C_{i-1})-p(C_{i-1})$ 
for $X\subseteq S_{i}$.  
In other words, $p_{i}$ is a set-function on $S_{i}$ 
obtained from $p$ by deleting $C_{i-1}$ and contracting $S-C_{i}$.  
The unique supermodular function $p_{F}$ 
defining the face $F$ is given by 
$\sum [p_{i}(S_{i}\cap X) :  i=1,\dots ,\ell]$.
A face $F$ is the set of elements $x$ of $B$ minimizing $cx$
whenever $c:S\rightarrow {\bf R} $ is a linear cost function such that
$c(s) = c(t)$ if $s,t\in S_{i}$ for some $i$ and $c(s) > c(t)$ 
if $s\in S_{i}$ and $t\in S_{j}$ for some subscripts $i<j$.

The intersection of an integral base-polyhedron 
$B=B'(p)$ \ ($=B(\overline{p})$) 
and an integral box $T(f,g)$ is an integral base-polyhedron.  
The intersection is non-empty 
if and only if
\begin{equation} 
p \leq \widetilde g \quad \hbox{and} \quad \widetilde f \leq \overline{p} .
\label{(pgfp)} 
\end{equation}

For an element $m$ of a base-polyhedron $B=B(b)$ defined 
by a submodular function $b$, we call a subset $X\subseteq S$ 
{\bf $m$-tight} (with respect to $b$) if $\widetilde m(X)=b(X)$.  
It is known (e.g.,  Theorem 14.2.8 in \cite{Frank-book}) 
that, for a given subset $X\subseteq S$, 
the face $B_{X} := \{ x\in B:  \widetilde x(X)=b(X) \}$ of $B$ 
is a non-empty base-polyhedron.  
This means that a subset $X$ is $m$-tight precisely if $m$ is in $B_{X}$.  
Clearly, the empty set and $S$ are $m$-tight, and $m$-tight sets are closed
under taking union and intersection
(see, for example, Lemma 2.4.7 in \cite{Frank-book}).
Therefore, for each subset $Z\subseteq S$, 
the intersection $T_{m}(Z;b)$ of all $m$-tight sets including $Z$ is 
the unique smallest $m$-tight set including $Z$. 
When $Z=\{s\}$ is a singleton, we simply write
$T_{m}(s;b)$ to denote the smallest $m$-tight set containing $s$.  
This set admits a representation
$T_{m}(s; b) = \{ t \in S: m+\chi_{s}-\chi_{t} \in B(b) \}$.
When the submodular function $b$ is 
understood from the context, 
we abbreviate $T_{m}(Z;b)$ to $T_{m}(Z)$.

Analogously, when $B=B'(p)$ is given by a supermodular function $p$,
we call $X\subseteq S$ {\bf $m$-tight} (with respect to $p$) if
$\widetilde m(X)=p(X)$.  
In this case, we also use the analogous notation 
$T_{m}(Z)=T_{m}(Z;p)$ and $T_{m}(s)=T_{m}(s;p)$.  
We have $T_{m}(s; p) = \{ t \in S: m-\chi_{s}+\chi_{t} \in B'(p) \}$.
Observe that for complementary functions $b$ and $p$, 
$X$ is $m$-tight with respect to
$b$ precisely if $S-X$ is $m$-tight with respect to $p$.

\begin{example} \rm \label{EXori3}
The set of in-degree vectors of orientations of a given undirected graph
forms an M-convex set.
Consider the undirected graph $G=(V,E)$ 
in the left-most of Figure \ref{FGoritri},
where $V=\{a,b,c \}$ and the set $E$ of edges 
consists of $ab$, $ac$, and a pair of parallel edges between $b$ and $c$.
An orientation of $G$ means 
a directed graph $D$ that is obtained from $G$ by orienting the edges of $G$.
The in-degree vector of $D$ is 
the vector $m$ on $V$ whose component at $v \in V$
is equal to the number of edges entering $v$ in $D$.
Three different orientations are depicted in Figure \ref{FGoritri}
with the corresponding in-degree vectors,
$m_{1}=(2,1,1)$,
$m_{2}=(2,2,0)$, and
$m_{3}=(1,3,0)$.
The set of in-degree vectors of all orientations of $G$
is known to form an M-convex set, say, $\odotZ{B}$.
In the present example, the M-convex set $\odotZ{B}$ consists of 10 members as
\[
\odotZ{B} = \{ \underline{(2,1,1)},  \underline{(2,2,0)},(2,0,2), 
\underline{(1,3,0)}, (1,0,3), (1,2,1), (1,1,2),
 (0,2,2), (0,3,1),  (0,1,3)    \},
\]
where  $m_{1}$, $m_{2}$, and $m_{3}$  are underlined.
Note that different orientations may result in the same in-degree vector.
The supermodular function $p$ describing the M-convex set $\odotZ{B}$
(or the base-polyhedron $B$)
is given as follows.
For any $X \subseteq V$, let $i_G(X)$ denote the number of edges of $G$ induced by $X$,
that is, 
$i_G(X) = | \{  uv \in E : \{ u, v \} \subseteq X \} |$. 
This function $i_{G}$ is a (non-negative) integer-valued supermodular function,
and we have $B=B'(p)$ for $p=i_{G}$. 
For any orientation $D$ with in-degree vector $m$, 
a subset $X$ of $V$ is $m$-tight
if and only if there is no edge entering $X$ in $D$.
In the present example,
$X = \{ b,c \}$ is $m$-tight for $m=m_{1}, m_{2}$,
whereas it is not $m_{3}$-tight.
\finbox
\end{example}


\section{Characterizing a decreasingly minimal element}
\label{SCchardecmin}

\subsection{Decreasing minimality}
\label{SCdefdecmin}

For a vector $x$, let $x {\downarrow}$
denote the vector obtained from $x$ by rearranging
its components in a decreasing order.  For example,
we call two vectors 
$x$ and $y$ (of same dimension) 
{\bf value-equivalent} if $x{\downarrow}= y{\downarrow}$.

A vector $x$ is {\bf decreasingly smaller} than vector $y$, in
notation $x <_{\rm dec} y$ \ if $x{\downarrow}$ \ is lexicographically
smaller than \ $y{\downarrow}$ in the sense that they are not
value-equivalent and $x{\downarrow}(j) < y{\downarrow}(j)$ for the
smallest subscript $j$ for which $x{\downarrow}(j)$ and
$y{\downarrow}(j)$ differ.  
For example, $x = (2,5,5,1,4)$ is decreasingly smaller than $y =(1,5,5,5,1)$
\ since $x{\downarrow}= (5,5,4,2,1)$  
is \ lexicographically smaller than \ $y{\downarrow}=(5,5,5,1,1)$.
We write $x\leq_{\rm dec} y$ to mean that $x$ is decreasingly smaller than or
value-equivalent to $y$.

For a set \ $Q$ \ of vectors,
$x\in Q$ is {\bf decreasingly minimal} 
({\bf dec-min}, for short) 
if \ $x \leq_{\rm dec} y$ \ for every \ $y\in Q$.  
Note that the dec-min elements of $Q$ are value-equivalent.  
Therefore an element $m$ of $Q$ is dec-min if
its largest component is as small as possible, within this, its second
largest component (with the same or smaller value than the largest
one) is as small as possible, and so on.  
An element $x$ of $Q$ is
said to be a {\bf max-minimized} element (a {\bf max-minimizer}, for short) 
if its largest component is as small as possible.  
A max-minimizer element $x$ is {\bf pre-decreasingly minimal}
 ({\bf pre-dec-min}, for short) in $Q$ 
if the number of its largest components is as small as possible.  
Obviously, a dec-min element is pre-dec-min, 
and a pre-dec-min element is max-minimized.
In Example~\ref{EXori3}, for example,
$m_{1}=(2,1,1)$ is dec-min in $\odotZ{B}$,
$m_{2}=(2,2,0)$ is a max-minimizer that is not dec-min,
and $m_{3}=(1,3,0)$ is not a max-minimizer.

In an analogous way, for a vector $x$, we let $x {\uparrow}$
denote the vector obtained from $x$ by rearranging its
components in an increasing order.  
A vector $y$ is {\bf increasingly larger} than vector $x$, 
in notation $y >_{\rm inc} x$, 
if they are not value-equivalent and 
$y{\uparrow}(j)> x{\uparrow}(j)$ 
holds for the smallest subscript $j$ for which 
$y{\uparrow}(j)$ and $x{\uparrow}(j)$ differ.  
We write $y \geq_{\rm inc} x$ 
if either $y >_{\rm inc} x$ or $x$ and $y$ are value-equivalent.  
Furthermore, we call an element $m$ of $Q$ 
{\bf increasingly maximal}
 ({\bf inc-max} for short) 
if its smallest component is as large as possible
over the elements of $Q$, within this its second smallest component is
as large as possible, and so on.

It should be emphasized that a dec-min element of a base-polyhedron
$B$ is not necessarily integer-valued.  
For example, if
$B=\{(x_{1},x_{2}):  \ x_{1}+x_{2}= 1\}$, 
then $x\sp{*}=(1/2,1/2)$ is a dec-min element of $B$.  
In this case, the dec-min members of $\odotZ{B}$ are
$(0,1) $ and $(1,0)$.

Therefore, finding a dec-min element of $B$ and finding a dec-min
element of $\odotZ{B}$ (the set of integral points of $B$) are two
distinct problems, and we shall concentrate only on the second,
discrete problem.  
In what follows, the slightly sloppy term integral
dec-min element of $B$ will always mean a dec-min element of $\odotZ{B}$.  
(The term is sloppy in the sense that an integral dec-min
element of $B$ is not necessarily a dec-min element of $B$).

We call an integral vector $x\in {\bf Z}\sp S$ {\bf uniform} 
if all of its components are the same integer $\ell$, 
and {\bf near-uniform} if
its largest and smallest components differ by at most 1, that is, if
$x(s)\in \{\ell, \ell +1\}$ for some integer $\ell$ for every $s\in S$.  
Note that if $Q$ consists of integral vectors and the
component-sum is the same for each member of $Q$, then any
near-uniform member
of $Q$ is obviously both decreasingly
minimal and increasingly maximal integral vector.

\subsection{Characterizing a dec-min element} 
\label{SCdecmincond}

Let $B=B(b)=B'(p)$ be a base-polyhedron defined by an integer-valued
submodular function $b$ or supermodular function $p$ 
(where $b$ and $p$ are complementary set-functions).  
Let $m$ be an integral member of $B$, that is, $m \in \odotZ{B}$.  
Recall the definition of $m$-tightness introduced at the end of Section~\ref{egal2base}.

The equivalences in the next claim will be used throughout.

\begin{claim}   \label{CLst}
Let $m \in \odotZ{B}$, and let
$s$ and $t$ be elements of $S$, and $m':=m+\chi_{s}-\chi_{t}$.
The following properties are pairwise equivalent.  
\smallskip

\noindent 
{\rm (A)} \ $m'\in \odotZ{B}$.  
\smallskip

\noindent 
{\rm (P1)} \ There is no $t\overline{s}$-set which is $m$-tight
with respect to $p$. 
\smallskip

\noindent 
{\rm (P2)} \ $s\in T_{m}(t;p)$.  
\smallskip

\noindent 
{\rm (B1)} \ There is no $s\overline{t}$-set which is $m$-tight
with respect to $b$.  
\smallskip

\noindent 
{\rm (B2)} \ $t\in T_{m}(s;b)$.  
\end{claim}

\Proof
(A) $\Rightarrow$ (P1) \ 
If $X$ is a $t\overline{s}$-set which is $m$-tight with respect to $p$, 
then $\widetilde m'(X)< \widetilde m(X) = p(X)$, 
showing that $m'\not \in \odotZ{B}$.

\smallskip

(P1) $\Rightarrow$ (P2) \ 
Since $T_{m}(t;p)$ is an $m$-tight set
containing $t$, (P1) implies that $s$ cannot be outside $T_{m}(t;p)$.

\smallskip

(P2) $\Rightarrow$ (A) \ 
Suppose that $m'$ is not in $\odotZ{B}$, that is, 
there is a subset $X\subset S$ with $\widetilde m'(X) < p(X)$.  
Since $m$ and $p$ are integer-valued and 
$\widetilde m(X) \geq p(X)$, 
we get from the definition of $m'$ that 
$\widetilde m(X)=p(X)$ and $X$ is a $t\overline{s}$-set, contradicting (P2).

\smallskip

(P1) $\Leftrightarrow $ (B1) \ 
Since $p$ and $b$ are complementary set-functions, 
a subset $X$ is $m$-tight with respect to $p$ precisely
if $S-X$ is $m$-tight with respect to $b$.

\smallskip

(P2) $\Leftrightarrow $ (B2) \ 
Suppose that $t\not \in T_{m}(s;b)$, that is, 
there is an $s\overline{t}$-subset $X$ which is $m$-tight with respect to $b$.  
Then $S-X$ is $m$-tight with respect to $p$, implying that
$s\not \in T_{m}(t;p)$, 
that is, (P2) implies (B2).  
The reverse implication follows analogously.  
\finbox

\medskip

A {\bf 1-tightening step} for $m\in \odotZ{B}$ 
is an operation that replaces $m$ by  
$m':=m+\chi_{s}-\chi_{t}$
where $s$ and $t$ are elements of $S$ for which $m(t)\geq m(s)+2$ and
$m'$ belongs to $\odotZ{B}$.  
Note that $m'$ is both decreasingly
smaller and increasingly larger than $m$.
Intuitively, a 1-tightening step may be viewed as a local improvement at $m$.
Since the mean of the components of $m$ does not change 
at a 1-tightening step while the square-sum of the components of $m$
strictly drops, consecutive 1-tightening steps may occur 
only a finite number of times (even if $B$ is unbounded).

As an example, consider the in-degree vector
$m=m_{3}=(1,3,0)$ in Example~\ref{EXori3},
where $S = \{ a,b,c \}$ and $m = (m(a), m(b), m(c))$.
For $(s,t)=(a,b)$ we have $m(t)\geq m(s)+2$ and 
$m'=m+\chi_{s}-\chi_{t} = (2,2,0) \in \odotZ{B}$.
(Note that $m'$ is equal to $m_{2}$ in Fig.~\ref{FGoritri}.)
Therefore, this is a 1-tightening step.
In contrast, for $(s,t)=(c,b)$ and $m=(1,3,0)$, we do not have 
a 1-tightening step since $m(t) =m(s)+1$, although 
$m+\chi_{s}-\chi_{t} = (0,3,1) \in \odotZ{B}$.

The next claim shows equivalent conditions
for the non-existence of a 1-tightening step.

\begin{claim}   \label{local} 
For an integral element $m$ of the integral base-polyhedron $B=B(b)=B'(p)$, 
the following conditions are pairwise equivalent.  
\smallskip

\noindent
{\rm (A)} 
\ There is no 1-tightening step for $m$.
\smallskip

\noindent 
{\rm (P1)} \ $m(s) \geq m(t)-1$ holds whenever $t\in S$ and
$s\in T_{m}(t;p)$.  
\smallskip

\noindent {\rm (P2)} \ 
Whenever $m(t)\geq m(s)+2$, there is a $t\overline{s}$-set $X$
which is $m$-tight with respect to $p$.  
\smallskip

\noindent 
{\rm (B1)} \ $m(s) \geq m(t)-1$ holds whenever $s\in S$ and $t\in T_{m}(s;b)$.
\smallskip

\noindent 
{\rm (B2)} \ 
Whenever $m(t)\geq m(s)+2$, there is an $s\overline{t}$-set $Y$
which is $m$-tight with respect to $b$.  
\finbox 
\end{claim}

\Proof
(A) $\Rightarrow$ (P1) \ 
If we had a pair $(s,t)$ of elements with $t\in S$, 
$s\in T_{m}(t;p)$ for which $m(s)\leq m(t)-2$, then replacing
$m$ by $m':= m + \chi_s - \chi_t$ would be a 1-tightening step.

\smallskip

(P1) $\Rightarrow$ (P2) \ 
Suppose that $(s,t)$ is a pair of elements 
for which $m(t)\geq m(s)+2$ but no $m$-tight $t\overline{s}$-set exists.  
Then $s\in T_{m}(t;p)$, contradicting (P1).

\smallskip

(P2) $\Rightarrow$ (A) \ 
Suppose there is a 1-tightening step for $m$, 
that is, there are elements $s$ and $t$ for which 
$m(t) \geq m(s) + 2$ and $m':=m+\chi_s-\chi_t$ is in $\odotZ{B}$.  
But $m'\in \odotZ{B}$ implies that 
no $m$-tight $t\overline{s}$-set can exist,
contradicting (P2).

\smallskip

(P1) $\Leftrightarrow $ (B1) \ follows from the equivalence of $s\in
T_{m}(t;p)$ and $t\in T_{m}(t;b)$ established in Claim~\ref{CLst}.

\smallskip

(P2) $\Leftrightarrow $ (B2) \ follows from the property that a subset
of $S$ is $m$-tight with respect to $b$ precisely if its complement is
$m$-tight with respect to $p$.  
\finbox

\medskip

For a given vector $m$ in ${\bf R}\sp S$, we call a set
$X\subseteq S$ an {\bf $m$-top set} (or a top-set with respect to $m$)
if $m(u)\geq m(v)$ holds whenever $u\in X$ and $v\in S-X$.  
Both the empty set and the ground-set $S$ are $m$-top sets, 
and $m$-top sets are closed under taking union and intersection.  
If $m(u)>m(v)$ holds
whenever $u\in X$ and $v\in S-X$, we speak of a {\bf strict $m$-top set}.  
A set $X \subseteq S$is a strict $m$-top set if and only if
$X$ is represented as $X = \{ s \in S : m(s) \geq \alpha \}$ for some integer $\alpha$.
For example, the vector $m = (4,2,2,1,1)$, indexed by 
$S=\{ s_{1},s_{2},s_{3},s_{4},s_{5} \}$,
has four strict $m$-top sets:
the empty set, 
$\{ s_{1} \}$, $\{ s_{1},s_{2},s_{3} \}$, and $S$.
Note that the number of strict non-empty $m$-top sets is at
most $n$ for every $m\in \odotZ{B}$ while $m \equiv 0$ exemplifies that
even all of the non-empty subsets of $S$ can be $m$-top sets.

\begin{theorem}  \label{equi.1} 
Let $b$ be an integer-valued submodular function 
and let $p:= \overline{b}$ be its complementary  (supermodular) function.  
For an integral element $m$ of the integral base-polyhedron $B=B(b)=B'(p)$,
the following four conditions are pairwise equivalent.
\smallskip

\noindent
{\rm (A)} \ There is no 1-tightening step for $m$
(or any one of the four
other equivalent properties holds in 
Claim {\rm \ref{local}}).

\smallskip

\noindent 
{\rm (B)} \ There is a chain $\cal C$ of $m$-top sets 
$(\emptyset \subset ) \ C_{1}\subset C_{2}\subset \cdots \subset C_\ell = S$ 
which are $m$-tight with respect to $p$ 
(or equivalently, whose complements are $m$-tight with respect to $b$) 
such that the restriction $m_{i}=m\vert S_{i}$ of $m$ to $S_{i}$ is
near-uniform for each member $S_{i}$ of the $S$-partition 
$\{S_{1},\dots,S_\ell\}$, where $S_{1}=C_{1}$ 
and $S_{i}:=C_{i}-C_{i-1}$ $(i=2,\dots ,\ell)$.
\smallskip

\noindent 
{\rm (C1)} \ $m$ is decreasingly minimal in $\odotZ{B}$.
\smallskip

\noindent 
{\rm (C2)} \ $m$ is increasingly maximal in $\odotZ{B}$.  
\end{theorem}

\Proof \ 
(B)$\rightarrow $(A): \ If $m(t)\geq m(s)+2$, then there is an
$m$-tight set $C_{i}$ containing $t$ and not containing $s$, 
from which Property (A) follows from Claim \ref{local}.

\smallskip

(A)$\rightarrow $(B): \ 
 Let $\cal C$ be a
longest chain consisting of non-empty $m$-tight and $m$-top sets
$C_{1}\subset C_{2}\subset \cdots \subset C_\ell =S$.  
For notational convenience, let 
$C_{0}=\emptyset $ (but $C_{0}$ is not a member of $\cal C$).  
We claim that $\cal C$ meets the requirement of (B).  
If, indirectly, this is not the case, then there is a subscript 
$i\in \{1,\dots ,\ell\}$ 
for which $m$ is not near-uniform within
$S_{i}:=C_{i}-C_{i-1}$.  
This means that the max $m$-value $\beta_{i}$ in
$S_{i}$ is at least 2 larger than the min $m$-value 
$\alpha_{i}$ in $S_{i}$, that is, 
$\beta_{i}\geq \alpha_{i} +2$.  
Let 
$Z:= \cup [ T_{m}(t;p): t\in S_{i}, \ m(t)=\beta_{i}]$.  
Then $Z$ is $m$-tight.  
Since $C_{i}$ is $m$-tight, $T_{m}(t;p)\subseteq C_{i}$ holds for $t\in S_{i}$ 
and hence $Z\subseteq C_{i}$.  
Furthermore, (A) implies that 
$m(v) \geq \beta_{i} -1$ for every $v\in Z\cap S_{i}$.

Consider the set $C':=C_{i-1}\cup Z$.  
Then $C'$ is $m$-tight, and
$C_{i-1}\subset C'\subset C_{i}$.  
Moreover, we claim that $C'$ is an $m$-top set.  
Indeed, if, indirectly, there is an element $u\in C'$ and an element $v\in S-C'$ 
for which $m(u)<m(v)$, 
then $u\in Z\cap S_{i}$ and 
$v\in C_{i}-Z$ since both $C_{i-1}$ and $C_{i}$ are $m$-top sets.  
But this is impossible since the $m$-value of each element of $Z\cap S_{i}$ 
is $\beta_{i}$ or $\beta_{i} -1$ while the $m$-value of each element
of $C_{i}-Z$ is at most $\beta_{i} -1$.

The existence of $C'$ contradicts the assumption that $\cal C$
was a longest chain of $m$-tight and $m$-top sets, and
therefore $m$ must be near-uniform within each $S_{i}$, that is, 
${\cal C}$ meets indeed the requirements in (B).

\smallskip

(C1)$\rightarrow $(A) and (C2)$\rightarrow $(A):  \ Property (A) must
indeed hold since a 1-tightening step for $m$ results in an element
$m'$ of $\odotZ{B}$ which is both decreasingly smaller and increasingly
larger than $m$.

\smallskip

(B)$\rightarrow $(C1): \ We may assume that the elements of $S$ are
arranged in an $m$-decreasing order 
$s_{1},\dots ,s_n$ 
(that is,
$m(s_{1})\geq m(s_{2})\geq \cdots \geq m(s_n)$) 
in such a way that each
$C_{i}$ in (B) is a starting segment.  
Let $m'$ be an element of
$\odotZ{B}$ which is decreasingly smaller than or value-equivalent to $m$.  
Recall that $m\vert X$ denoted the vector $m$ restricted to a subset $X\subseteq S$.

\begin{lemma}
For each $i=0,1,\dots ,\ell$, vector $m'\vert C_{i}$ is
value-equivalent to vector $m\vert C_{i}$.  
\end{lemma}  

\Proof Induction on $i$.  
For $i=0$, the statement is void so we assume that $1\leq i\leq \ell$.  
By induction, we may assume that the
statement holds for $j\leq i-1$ and we want to prove it for $i$.
Since $m'\vert C_{i-1}$ is value-equivalent to $m\vert C_{i-1}$ and
$C_{i-1}$ is $m$-tight, it follows that $C_{i-1}$ is $m'$-tight, too.

Let $\beta_{i}$ denote the max $m$-value of the elements of
$S_{i}=C_{i}-C_{i-1}$.  By the hypothesis in (B), 
the maximum and the minimum of the $m$-values in $S_{i}$ differ by at most 1. 
Hence we can assume that there are $r_{i}>0$ elements in $S_{i}$ 
with $m$-value $\beta_{i}$ and 
$\vert S_{i}\vert -r_{i}\geq 0$ elements 
with $m$-value $\beta_{i} -1$.

As $m\vert C_{i-1}$ is value-equivalent to $m'\vert C_{i-1}$ and $m'$
was assumed to be decreasingly smaller than or value-equivalent to
$m$, we can conclude that $m'\vert (S-C_{i-1})$ is decreasingly smaller 
than or value-equivalent to $m\vert (S-C_{i-1})$.  
Therefore, $S_{i}$ contains at most $r_{i}$ elements of $m'$-value $\beta_{i}$ and hence
\begin{align*}
p(C_{i}) &\leq \widetilde m'(C_{i}) = \widetilde m'(C_{i-1}) + \widetilde m'(S_{i}) 
\\ &
\leq \widetilde m'(C_{i-1}) + r_{i}\beta_{i} + (\vert S_{i}\vert -r_{i}) (\beta_{i} -1) 
\\ & 
=  \widetilde m(C_{i-1}) + r_{i}\beta_{i}
    + (\vert S_{i}\vert -r_{i}) (\beta_{i} -1) 
\\ & 
 = \widetilde m(C_{i-1})  + \widetilde m(S_{i}) = \widetilde m(C_{i}) 
 = p(C_{i}),
\end{align*}
from which equality follows everywhere.  
In particular, $S_{i}$ contains exactly $r_{i}$
elements of $m'$-value $\beta_{i}$ and $\vert S_{i}\vert -r_{i}$ elements
of $m'$-value $\beta_{i} -1$, proving the lemma.  
\finbox 
\medskip

By the lemma, $m'$ is value-equivalent to $m$, and hence $m$ is a
decreasingly minimal element of $\odotZ{B}$, that is, (C1) follows.

\smallskip

(B)$\rightarrow $(C2): \ The property in (C1) that $m$ is 
decreasingly minimal 
in $\odotZ{B}$ is equivalent to the statement that
$-m$ is increasingly maximal in $-{\odotZ{B}}$, 
that is, (C2) holds with respect to $-m$ and $-{\odotZ{B}}$.  
As we have already proved the implications 
(C2)$\rightarrow$(A)$\rightarrow$(B)$\rightarrow $(C1), 
it follows that (C1) holds for $-m$ and $-{\odotZ{B}}$.  
But (C1) for $-m$ and $-{\odotZ{B}}$ is just the same
as (C2) for $m$ and $\odotZ{B}$.  
\finbox \finboxHere
\medskip

\begin{figure}
\centering
\includegraphics[width=0.55\textwidth,clip]{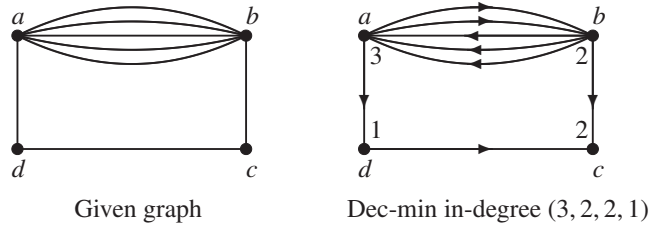}
\caption{Orientation of a graph (Example \ref{EXori4})} 
\label{FGorisqu}
\end{figure}

\begin{example} \rm \label{EXori4}
Theorem \ref{equi.1} is illustrated here for a graph orientation problem.
Consider the undirected graph $G=(V,E)$ 
in the left of Figure \ref{FGorisqu},
where $V=\{a,b,c,d\}$ and the set $E$ of edges 
consists of $bc$, $cd$, $da$, and five parallel edges between $a$ and $b$.
Recall that the associated M-convex set $\odotZ{B}$ consists of the in-degree vectors
of all orientations of $G$.
Let $D$ denote the orientation of $G$ depicted in the right of 
Figure \ref{FGorisqu}.
We shall apply  Theorem \ref{equi.1} to verify that 
the in-degree vector 
$m=(3,2,2,1)$ of $D$ is a dec-min element of $\odotZ{B}$.

To verify Condition (A) for $m$, first note that 
$m+\chi_s-\chi_t$ belongs to $\odotZ{B}$ for $s,t \in V$
precisely if 
there is a dipath from $s$ to $t$ in $D$.
For the in-degree vector $m=(3,2,2,1)$,
we have $m(t) \geq m(s) + 2$ 
only for $(s,t)=(d,a)$,
whereas there is no dipath from $d$ to $a$ in $D$.
Therefore, no 1-tightening step exists for $m=(3,2,2,1)$.

To verify Condition (B) for $m$, recall 
(from Example \ref{EXori3})
that a subset $X$ of $V$ is $m$-tight
if and only if there is no edge entering $X$ in $D$.
The $m$-tight sets are $\{ a,b \}$ and $\{ a,b,d \}$
as well as the empty set and $V$.
These $m$-tight sets are $m$-top sets except for $\{ a,b,d \}$.
Thus the (longest) chain of non-empty $m$-top and $m$-tight sets 
consists of $C_{1} \subset C_{2}$ where
$C_{1} = \{ a,b \}$ and $C_{2} = \{ a,b,c,d \}$.
This chain determines the partition of $V$ into two parts
$S_{1} = \{ a,b \}$ and $S_{2} = \{ c,d \}$,
for which 
$m_{1} = m | S_{1} = (3,2)$ and
$m_{2} = m | S_{2} = (2,1)$
are both near-uniform.
Therefore, Condition (B) is satisfied by $m=(3,2,2,1)$.

To verify Conditions (C1) and (C2), we may enumerate 
all possible in-degree vectors.
The subgraph consisting of three edges $bc$, $cd$ and $da$
admits six in-degree vectors:
\[
\odotZ{B_{0}}
 = \{ 
(0,1,1,1), (1,0,1,1), (1,1,1,0), (1,1,0,1), (0,0,2,1), (0,0,1,2) \}.
\]
By adding the five parallel edges we obtain
\[
\odotZ{B} = \{ m_{0} + (k,5-k,0,0)  : m_{0} \in \odotZ{B_{0}}, \ 0 \leq k \leq 5 \}.
\]
A straightforward inspection reveals that 
$m=(3,2,2,1)$ is dec-min and inc-max in $\odotZ{B}$.
Note that there are  three further dec-min elements: 
$(2,3,2,1)$, $(2,3,1,2)$, and $(3,2,1,2)$.
\finbox
\end{example}

\begin{remark} \rm \label{RMdemminincmax}
The equivalence of (C1) and (C2) in Theorem \ref{equi.1} 
shows that an element of an M-convex set is decreasingly minimal 
if and only if it is increasingly maximal.
In the intersection of two M-convex sets
(called an M$_{2}$-convex set in \cite{Murota03}),
however,
decreasing minimality  and increasing maximality do not coincide.
For example, consider two M-convex sets
\begin{align*} 
\odotZ{B_{1}} &= \{(2, 0, 0, 0), \ (1, -1, 1, 1), \ (2, -1, 1, 0), \ (1, 0, 0, 1)\}, 
\\ 
\odotZ{B_{2}} &= \{(2, 0, 0, 0), (1, -1, 1, 1), (2, \-1, 0, 1), (1, 0, 1, 0)\}.
\end{align*}
In their intersection 
$\odotZ{B_{1}} \cap \odotZ{B_{2}} = \{(2, 0, 0, 0), \ (1, -1, 1, 1)\}$,
the element $x = (2, 0, 0, 0)$ is increasingly maximal 
while $y = (1, -1, 1, 1)$ is decreasingly minimal.
\finbox
\end{remark}

\subsection{Minimizing the sum of the $k$ largest components}
\label{SCklargestsum}

A decreasingly minimal element of $\odotZ{B}$ has the starting property
that its largest component is as small as possible.  
As a natural extension, one may be interested in finding a member of $\odotZ{B}$ 
for which the sum of the $k$ largest components is as small as possible.
We refer to this problem as 
{\bf min $k$-largest-sum}.

\begin{theorem}\label{klargest} 
Let $B$ be an integral base-polyhedron and
$k$ an integer with $1\leq k\leq n$.  Then any dec-min element $m$ of
$\odotZ{B}$ is a solution to the min $k$-largest-sum problem.
\end{theorem}

\Proof 
Observe first that if $z_{1}$ and $z_{2}$ are dec-min elements of $\odotZ{B}$, 
then it follows from the very definition of decreasing
minimality that the sum of the first $j$ largest components of $z_{1}$
and of $z_{2}$ are the same for each $j=1,\dots ,n$.  

Let $K$ denote the
sum of the first $k$ largest components of any dec-min element, and
assume indirectly that there is a member $y\in \odotZ{B}$ for which the
sum of its first $k$ largest components is smaller than $K$.  
Assume that the componentwise square-sum of $y$ is as small as possible.  
By the previous observation, $y$ is not a dec-min element.  
Theorem~\ref{equi.1} implies that there are elements $s$ and $t$ of $S$ for
which $y(t)\geq y(s)+2$ and $y':=y - \chi_{t} + \chi_{s}$ is in $\odotZ{B}$.  
The sum of the first $k$ largest components of $y'$ is at
most the sum of the first $k$ largest components of $y$, and hence
this sum is also smaller than $K$.  
But this contradicts the choice of
$y$ since the componentwise square-sum of $y'$ is 
strictly smaller than that of $y$.  
\finbox 
\medskip

This theorem shows that M-convex sets have a striking property.  
Namely, any dec-min element of an M-convex set $\odotZ{B}$
is simultaneously a solution to the min $k$-largest-sum problem for
{\em each} $k=1,2,\dots ,n$.  
We say that such an element $m$ is a  
{\bf simultaneous $k$-largest-sum minimizer}.
This notion has been investigated 
in the literature of majorization \cite{AS18,MOA11,Tamir95}
under the name of `least majorized' element.
In particular, Tamir \cite{Tamir95}
proved the existence of a least majorized integral element for integral base-polyhedra. 
(Actually, he proved this even for g-polymatroids, but this more general result
is an easy consequence of the special case concerning base-polyhedra).

The following result shows that this property actually characterizes
dec-min elements of an M-convex set.

\begin{theorem}  \label{THmajordecmin} 
Let $B$ be an integral base-polyhedron.  An element $m$ of
$\odotZ B$ is dec-min if and only if $m$ is a simultaneous
$k$-largest-sum minimizer. 
\end{theorem}  

\Proof The content of Theorem \ref{klargest} is 
that a dec-min element is a simultaneous $k$-largest-sum minimizer.
To see the converse, let $m\in \odotZ{B}$ be a simultaneous
$k$-largest-sum minimizer.  
Suppose indirectly that $m$ is not dec-min.  
By Theorem \ref{equi.1}, there is a 1-tightening step for $m$, 
that is, 
there are elements $s$ and $t$ of $S$ with $m(s)\geq m(t)+2$ such
that $m':=m- \chi_{s} + \chi_{t}$ is in $\odotZ{B}$.  
Let $k'$ denote the number of components of $m$ with value at least $m(s)$.  
Then the sum of the $k'$ largest components of $m'$ is 
one less than the sum of the $k'$ largest components of $m$, 
contradicting the assumption that $m$ is $k$-sum-minimizer 
for each $k=1,2,\dots ,n$.  
\finbox 

\medskip

While the $k$-largest sum is a natural objective function
to conceive from the definition of dec-minimality, 
it is also possible to characterize the dec-min elements of an M-convex set
as minimizers of other convex functions,
which we shall discuss in depth in Section~\ref{SCsqdiffsum}.
We now turn to investigating the structure of 
the set of dec-min elements of an M-convex set.
In particular, our next goal is to introduce `dual' objects for dec-minimality,
which we call the canonical chain and the canonical partition. 
These `dual' objects are constructed 
in Sections \ref{peak} and \ref{canonical} below.

\section{Decomposition by pre-decreasingly minimal elements}
\label{peak}

We continue to assume that $p$ is an integer-valued (with possible
$-\infty $ values but with finite $p(S)$) 
supermodular function, which implies that $B=B'(p)$ is a non-empty
integral base-polyhedron.

One of the main goals of this paper 
is to show that the set ${\rm dm}(\odotZ{B})$ of all
dec-min elements of $\odotZ{B}$ is an M-convex set, meaning that
there exists an integral base-polyhedron 
$B\sp{\bullet} \subseteq B$ 
(obtained by intersecting a face of $B$ with a `small' integral box) 
such that ${\rm dm}(\odotZ{B})$ 
is the set of integral elements of $B\sp{\bullet}$. 
In addition, we shall show that 
$B\sp{\bullet}$ is a special base-polyhedron
which is obtained from a matroid base-polyhedron by translating it
with an integral vector.
The base-polyhedron $B\sp{\bullet}$ will be obtained with the help of a
decomposition of $B$ along a certain \lq canonical\rq\ partition
$\{S_{1},S_{2},\dots ,S_{q}\}$ of $S$ into non-empty sets.  
In this section, we construct 
the first member $S_{1}$ of this partition along
with a matroid on $S_{1}$.  
The set $S_{1}$, depending only on $B$, will
be called the peak-set of $B$.
The peak-set $S_{1}$ allows us to decompose the problem of
finding dec-min elements of $\odotZ{B}$ to two independent problems 
on $S_{1}$ and on $S-S_{1}$.
In Section~\ref{canonical}, we shall construct the other members 
$S_{2}, S_{3}, \dots$ of the partition 
by applying the same procedure iteratively.

\subsection{Max-minimizers and pre-dec-min elements}
\label{SCmaxminzerpredm}

In Section~\ref{SCdefdecmin} we introduced three related notions,
dec-min elements,
pre-dec-min elements, and 
max-minimizers. 
A dec-min element is always pre-dec-min, and 
a pre-dec-min element is always a max-minimizer.
Recall that an element of $\odotZ{B}$ is called a max-minimizer if 
its largest component is as small as possible, 
while a max-minimizer is  a pre-dec-min element of $\odotZ{B}$ 
if the number of its maximum components is smallest possible.

For a number  $\beta$, we say that a vector is {\bf $\beta$-covered} 
if each of its components is at most $\beta$.
Throughout our discussion, 
$\beta_{1}$ denotes the smallest integer for which
$\odotZ{B}$ has a $\beta_{1}$-covered  element, that is, 
\begin{equation} 
\beta_{1}:=  \min \{ \,  \max \{ z(s) : s \in S \} : z \in \odotZ{B} \, \} .
\label{(beta1-def)} 
\end{equation}
\noindent
Note that an element $m$ of $\odotZ{B}$ is $\beta_{1}$-covered 
precisely if $m$ is a max-minimizer. 
Moreover, $\beta_{1}$ is equal to the largest component of any 
pre-dec-min (and hence any dec-min) element of $\odotZ{B}$.  
For any real number $\alpha \in {\bf R}$, 
let $\lceil \alpha \rceil$ denote the smallest integer not smaller than $\alpha$.

\begin{theorem}   \label{betamin} 
For the largest component $\beta_{1}$ of a
max-minimizer of $\odotZ{B}$, one has 
\begin{equation} 
\beta_{1} =\max \{ \llceil {p(X)  \over \vert X\vert }\rrceil : 
    \emptyset \not =X\subseteq S\}.
\label{(betamin)} 
\end{equation} 
\end{theorem}

\Proof 
Formula \eqref{(pgfp)}, when applied to the special case 
with $f\equiv -\infty$ and $g\equiv \beta$, implies 
that $B$ has a $\beta$-covered element 
if and only if 
\begin{equation} 
\beta \vert X\vert \geq p(X) \quad
\hbox{whenever}\ \ X\subseteq S. 
\label{(beta-covered)} 
\end{equation} 
\noindent
Moreover,
if $\beta$ is an integer and \eqref{(beta-covered)} holds, then $B$
has an integral $\beta$-covered element.  
As $\beta \vert X\vert \geq
p(X)$ holds for an arbitrary $\beta$ when $X=\emptyset $, it follows
that the smallest integer $\beta$ meeting this \eqref{(beta-covered)}
is indeed 
$\max \{ \lceil {p(X) / \vert X\vert }\rceil : \emptyset \not =X\subseteq S\}$.  
\finbox 
\medskip

For a $\beta_{1}$-covered element $m$ of $\odotZ{B}$, let $r_{1}(m)$
denote the number of $\beta_{1}$-valued components of $m$.  Recall that
for an element $s\in S$ we denoted the unique smallest $m$-tight set
containing $s$ by $T_{m}(s)=T_{m}(s;p)$ (that is, $T_{m}(s)$ is the
intersection of all $m$-tight sets containing $s$).  
Furthermore, let
\begin{equation}  \label{(S1m.def)} 
S_{1}(m):=\cup \{T_{m}(t):  m(t)=\beta_{1}\}.  
\end{equation}
\noindent
Then $S_{1}(m)$ is $m$-tight and $S_{1}(m)$ is actually the unique
smallest $m$-tight set containing all the $\beta_{1}$-valued elements of $m$.

\begin{theorem}   \label{optkrit} 
A $\beta_{1}$-covered element $m$ of $\odotZ{B}$ is pre-dec-min 
if and only if 
$m(s)\geq \beta_{1} -1$ for each $s\in S_{1}(m)$.  
\end{theorem}

\Proof 
Necessity.  
Let $m$ be a pre-dec-min element of $\odotZ{B}$.
For any $\beta_{1}$-valued element $t\in S$ and any element $s\in
T_{m}(t)$, we claim that $m(s)\geq \beta_{1} -1$.  
Indeed, if we had $m(s)\leq \beta_{1} -2$, then the vector $m'$ arising from $m$ by
decreasing $m(t)$ by 1 and increasing $m(s)$ by 1 belongs to $B$
(since $T_{m}(t)$ is the smallest $m$-tight set containing $t$) 
and has one less $\beta_{1}$-valued components than $m$ has, contradicting the
assumption that $m$ is pre-dec-min.

Sufficiency.  
Let $m'$ be an arbitrary $\beta_{1}$-covered integral
element of $B$.  Abbreviate $S_{1}(m)$ by $Z$ and let $h'$ denote the
number of elements $z\in Z$ for which $m'(z)=\beta_{1}$.  
Then
\begin{eqnarray*} 
& &\vert Z\vert (\beta_{1} -1) + r_{1}(m) 
= \widetilde m(Z) = p(Z) \leq \widetilde m'(Z) 
\\ & & 
\ \leq h'\beta_{1} + (\vert Z\vert -h')(\beta_{1} -1) 
       = \vert Z\vert (\beta_{1} -1) + h' 
\\ & &
\ \leq  \vert Z\vert (\beta_{1} -1) + r_{1}(m') ,
\end{eqnarray*}
\noindent
from which $r_{1}(m)\leq r_{1}(m')$, as required.  
\finbox
\medskip

Define the set-function $h_{1}$ on $S$ as follows.  
\begin{equation} 
h_{1}(X):= p(X) - (\beta_{1} -1)\vert X\vert \ \ \hbox{for }\ \ X\subseteq S.
\label{(h1.def)} 
\end{equation}

\begin{theorem}   \label{betavalmin} 
For the minimum number $r_{1}$ of $\beta_{1}$-valued components 
of a $\beta_{1}$-covered member of $\odotZ{B}$,  one has
\begin{equation} 
r_{1}= \max \{h_{1}(X):  X\subseteq S\}.  
\label{(r1)} 
\end{equation} 
\end{theorem}

\Proof 
Let $m$ be an element of $\odotZ{B}$ for which the maximum of
its components is $\beta_{1}$, and let $X$ be an arbitrary subset of $S$.  
Suppose that $X$ has $\ell$ \ $\beta_{1}$-valued components.  
Then
\begin{equation} 
p(X) \ \leq  \ \widetilde m(X) 
   \ \leq \  \ell \beta_{1} + (\vert X\vert -\ell)(\beta_{1} -1) 
  \ = \  \vert X\vert (\beta_{1} -1) + \ell 
  \ \leq \  \vert X\vert (\beta_{1} -1) + r_{1}(m) ,
\label{(estim)} 
\end{equation}
\noindent
from which $r_{1}(m) \geq p(X) - (\beta_{1} -1)\vert X\vert =h_{1}(X)$, 
implying that 
\[
 r_{1} = \min \{r_{1}(m) :  m\in \odotZ{B},   
    \mbox{ \ $m$ is  $\beta_{1}$-covered } \} 
  \ \geq \  \max \{h_{1}(X): X\subseteq S\}. 
\]

In order to prove the reverse inequality, we have to find 
a $\beta_{1}$-covered integral element $m$ of $B$ and a subset $X$ of $S$ for
which $r_{1}(m) = h_{1}(X)$, 
which is equivalent to requiring that each of the three
inequalities in \eqref{(estim)} holds with equality.  
That is, the following three {\bf optimality criteria} hold:  
\ (a) \ $X$ is $m$-tight, 
\ (b) \ $X$ contains all $\beta_{1}$-valued components of $m$, and 
\ (c) \ $m(s)\geq \beta_{1} -1$ for each $s\in X$.

Let $m$ be a pre-dec-min element of $B$.  
Then $S_{1}(m)$ is $m$-tight, $S_{1}(m)$ contains all $\beta_{1}$-valued elements and, 
by Theorem~\ref{optkrit}, $m(s)\geq \beta_{1} -1$ for all $s\in S_{1}(m)$, 
therefore $m$ and $S_{1}(m)$ satisfy the three optimality criteria.  
\finbox

\medskip

Note that $r_{1}$ is the number of $\beta_{1}$-valued components of any
pre-dec-min element (and in particular, any dec-min element) of $\odotZ{B}$.

\subsection{The peak-set $S_{1}$}

Since the set-function $h_{1}$ introduced in \eqref{(h1.def)} is supermodular, 
the maximizers of $h_{1}$ are closed under taking intersection and union.  
Let $S_{1}$ denote the unique smallest subset of $S$ maximizing $h_{1}$.  
In other words, $S_{1}$ is the intersection of all sets maximizing $h_{1}$.  
We call this set $S_{1}$ the {\bf peak-set} of $B$ (and of $\odotZ{B}$).

\begin{theorem}  \label{smallest} 
For every pre-dec-min 
(and in particular, for every dec-min) 
element $m$ of $\odotZ{B}$, the set $S_{1}(m)$
introduced in \eqref{(S1m.def)} is independent of the choice of $m$
and $S_{1}(m)=S_{1}$, where $S_{1}$ is the peak-set of $B$.  
\end{theorem}

\Proof 
It follows from Theorem~\ref{betavalmin} that, given a
pre-dec-min element $m$ of $B$, a subset $X$ is maximizing $h_{1}$ 
precisely if the three optimality criteria mentioned in the proof hold.  
Since $S_{1}(m)$ meets the optimality criteria, 
it follows that $S_{1}\subseteq S_{1}(m)$.  
If, indirectly, there is an element $s\in S_{1}(m) -S_{1}$, 
then $m(s)=\beta_{1} -1$ since $S_{1}$ contains all the
$\beta_{1}$-valued elements.  
By the definition of $S_{1}(m)$, there is a $\beta_{1}$-valued element 
$t\in S_{1}(m)$ for which the smallest $m$-tight set $T_{m}(t)$ contains $s$, 
but this is impossible since $S_{1}$ is an $m$-tight set containing $t$ but not $s$.  
\finbox

\medskip

Since $S_{1}=S_{1}(m)$ is $m$-tight and near-uniform, we obtain that 
\[
\beta_{1}
= \llceil { \widetilde m (S_{1}) \over \vert S_{1}\vert }\rrceil 
= \llceil {p(S_{1}) \over \vert S_{1}\vert }\rrceil , 
\]
\noindent
and the definitions of $S_{1}$ and $r_{1}$ imply that 
\begin{equation} 
  r_{1}=p(S_{1}) - (\beta_{1}-1)\vert S_{1}\vert . 
\label{(r11)} 
\end{equation}

\begin{proposition}   \label{PR=S1} 
$S_{1}=\{s\in S:  \mbox{\rm  there is a pre-dec-min element } \
       m \in \odotZ{B} \ \mbox{\rm with } \  m(s)=\beta_{1}\}$.  
For every pre-dec-min element $m$ of $\odotZ{B}$,
$m(s)\geq \beta_{1} -1$ holds for every $s\in S_{1}$, and  
$m(s)\leq \beta_{1} -1$ holds for every $s\in S-S_{1}$.
\end{proposition}

\Proof 
If $m(s)=\beta_{1}$ for some pre-dec-min $m$, then $s\in S_{1}(m)=S_{1}$.  
Conversely, let $s\in S_{1}$ and let $m$ be a pre-dec-min element.  
We are done if $m(s)=\beta_{1}$.  
If this is not the case,
then $m(s)=\beta_{1} -1$ by Theorem~\ref{optkrit}.  
By the definition of $S_{1}(m)$, 
there is an element $t\in S_{1}(m)$ for which $m(t)=\beta_{1}$ and $s\in T_{m}(t)$.  
But then $m':= m+\chi_{s}-\chi_{t}$
 is in $\odotZ{B}$, $m'(s)=\beta_{1}$ and $m'$ is also
pre-dec-min as it is value-equivalent to $m$.  
\finbox

\subsection{Separating along $S_{1}$}
\label{SCsepS1}

Let $S_{1}$ be the peak-set occurring in Theorem~\ref{smallest} and 
let $S_{1}':=S-S_{1}$.  
Let $p_{1}=p\vert S_{1}$ denote the restriction of $p$ to
$S_{1}$, and let $B_{1}\subseteq {\bf R}\sp {S_{1}}$ denote the
base-polyhedron defined by $p_{1}$, that is, $B_{1}:=B'(p_{1})$.  
Suppose that $S_{1}'\not =\emptyset $ and let $p_{1}':=p/S_{1}$, 
that is, $p_{1}'$ is the set-function on $S_{1}'$ 
obtained from $p$ by contracting $S_{1}$ \
($p_{1}'(X)= p(S_{1}\cup X)-p(S_{1})$ for $X\subseteq S_{1}'$).

Consider the face $F$ of $B$ determined by $S_{1}$, that is, 
$F$ is the direct sum of the base-polyhedra $B_{1}=B'(p_{1})$ and $B_{1}'=B'(p_{1}')$.
Then the dec-min elements of $\odotZ{B}_{1}$ are exactly the integral
elements of the intersection of $B_{1}$ and the box given by 
$\{x: \beta_{1} -1\leq x(s)\leq \beta_{1} \ \mbox{ for every } s\}$.  
 Hence the dec-min elements of
$\odotZ{B}_{1}$ are near-uniform.

\begin{theorem}   \label{felbont} 
An integral vector $m=(m_{1},m_{1}')$ is a dec-min element of $\odotZ{B}$ 
if and only if 
$m_{1}$ is a dec-min element of $\odotZ{B}_{1}$ 
and $m_{1}'$ is a dec-min element of $\odotZ{B_{1}'}$.  
\end{theorem}

\Proof 
Suppose first that $m$ is a dec-min element of $\odotZ{B}$.
Then $S_{1}=S_{1}(m)$ by Theorem~\ref{smallest} and $m$ is a max-minimizer,
implying that every component of $m$ in $S_{1}$ is of value $\beta_{1} -1$
or value $\beta_{1}$, and $m$ has exactly $r_{1}$ components of value $\beta_{1}$.  
Therefore each of the components of $m_{1}$ is $\beta_{1} -1$ or $\beta_{1}$, 
that is, $m_{1}$ is near-uniform.
Since $m_{1}$ is obviously in $\odotZ{B}_{1}$,
$m_{1}$ is indeed dec-min in $\odotZ{B}_{1}$.

Since $\widetilde m(S_{1}) = p(S_{1})$, for a set $X\subseteq S_{1}'$, we have 
\[
\widetilde m_{1}'(X) = \widetilde m(X) 
  = \widetilde m(S_{1}\cup X)- \widetilde m(S_{1}) 
  = \widetilde m(S_{1}\cup X) - p(S_{1}) 
  \geq p(S_{1}\cup X) - p(S_{1}) = p_{1}'(X).
\]
\noindent
Furthermore 
\[
\widetilde m_{1}'(S_{1}') 
= \widetilde m(S_{1}') 
= \widetilde m(S_{1}\cup S_{1}') - \widetilde m(S_{1}) 
= p(S_{1}\cup S_{1}') - p(S_{1}) = p_{1}'(S_{1}'),
\] 
that is, $m_{1}'$ is in $\odotZ{B_{1}'}$.  If,
indirectly, $m_{1}'$ is not dec-min, then, by applying Theorem
\ref{equi.1} to $S_{1}'$, $m_{1}'$, and $p_{1}'$, we obtain that there are
elements $t$ and $s$ of $S_{1}'$ for which $m_{1}'(t) \geq m_{1}'(s)+2$ and
\ $(*)$ \ no $t\overline{s}$-set exists which is $m_{1}'$-tight with respect to $p_{1}'$.  
On the other hand, $m$ is a dec-min element of $\odotZ{B}$ for which 
\[
m(t)=m_{1}'(t) \geq m_{1}'(s)+2 = m(s)+2,
\] 
and hence there must
be a $t\overline{s}$-set $Y$ which is $m$-tight with respect to $p$.

Since $S_{1}$ is $m$-tight with respect to $p$, 
the set $S_{1}\cup Y$ is also $m$-tight with respect to $p$.  
Let $X:=S_{1}'\cap Y$.  
Then
\[
\widetilde m(X) + \widetilde m(S_{1}) 
  = \widetilde m(S_{1}\cup Y) 
  = p(S_{1}\cup Y) = p (S_{1}\cup X),
\]
and hence
\[
\widetilde m_{1}'(X)  = \widetilde m(X) 
= p(S_{1}\cup X) - \widetilde m(S_{1}) = p(S_{1}\cup X) - p(S_{1}) = p_{1}'(X),
\]
that is, $X$ is a $t\overline{s}$-set which is $m_{1}'$-tight with respect to $p_{1}'$, 
in contradiction with statement $(*)$ above that no such set exists.

To see the converse, assume that $m_{1}$ is a dec-min element of
$\odotZ{B}_{1}$ and $m_{1}'$ is a dec-min element of $\odotZ{B_{1}'}$.  
This immediately implies that $m$ is in the face $F$ of $B$ determined by $S_{1}$.  
Suppose, indirectly, that $m$ is not a dec-min element of $\odotZ{B}$.  
By Theorem~\ref{equi.1}, there are elements $t$ and $s$
of $S$ for which $m(t) \geq m(s)+2$ and 
$(**)$ no $t\overline{s}$-set exists which is $m$-tight with respect to $p$.  
If $t\in S_{1}$, then $s$
cannot be in $S_{1}$ since the $m$-value of each element of $S_{1}$ is
$\beta_{1}$ or $\beta_{1} -1$.  But $S_{1}$ is $m_{1}$-tight with respect to
$p$ and hence it is $m$-tight with respect to $p$, contradicting property $(**)$.  
Therefore $t$ must be in $S_{1}'$, 
implying, by Proposition \ref{PR=S1}, that $s$ is also in $S_{1}'$.

Since $m_{1}'$ is a dec-min element of $\odotZ{B_{1}'}$, there must be a
$t\overline{s}$-set $Y\subset S_{1}'$ which is $m_{1}'$-tight with respect to $p_{1}'$.  
It follows that 
\[
\widetilde m(Y) = \widetilde m_{1}'(Y) 
= p_{1}'(Y) = p(S_{1}\cup Y) - p(S_{1}) 
\leq \widetilde m(S_{1}\cup Y) - \widetilde m(S_{1})=\widetilde m(Y),
\]
from which $\widetilde m(S_{1}\cup Y) = p(S_{1}\cup Y)$, 
contradicting property $(**)$ that no $t\overline{s}$-set exists which is $m$-tight with respect to $p$.  
\finbox 
\medskip

An important consequence of Theorem~\ref{felbont} is that, in order to
find a dec-min element of $\odotZ{B}$, it will suffice to find
separately a dec-min element of $\odotZ{B_{1}}$ (which was shown above to
be a near-uniform vector) and a dec-min element of $\odotZ{B_{1}'}$.  
The algorithmic details are discussed in \cite{FM19partB}.

\begin{theorem}  \label{m1.ekv} 
Let $S_{1}$ be the peak-set of $\odotZ{B}$. 
For an element $m_{1}$ of $\odotZ{B_{1}}$, 
the following properties are pairwise equivalent.
\smallskip

\noindent 
{\rm (A1)} \ 
$m_{1}$ has $r_{1} \ \ (= p(S_{1}) - (\beta_{1}-1)\vert S_{1}\vert $ \ $>0)$\ \ 
components of value $\beta_{1}$ 
and 
$\vert S_{1}\vert -r_{1}$ \ $(\geq 0)$ \ 
components of value $\beta_{1} -1$.

\smallskip

\noindent 
{\rm (A2)} \ $m_{1}$ is near-uniform. 

\smallskip

\noindent 
{\rm (A3)} \ $m_{1}$ is dec-min in $\odotZ{B_{1}}$.  
\smallskip

\noindent 
{\rm (B1)} \ $m_{1}$ is the restriction of a dec-min element
$m$ of $\odotZ{B}$ to $S_{1}$.  
\smallskip

\noindent 
{\rm (B2)} \ $m_{1}$ is the restriction of a pre-dec-min
element $m$ of $\odotZ{B}$ to $S_{1}$.  
\end{theorem}

\Proof 
The implications (A1)$\rightarrow $(A2)$\rightarrow $(A3) and
(B1)$\rightarrow $(B2) are immediate from the definitions.

(A3)$\rightarrow $(B1):  
Let $m_{1}'$ be an arbitrary dec-min element of $\odotZ{B_{1}'}$.  
By Theorem~\ref{felbont}, $m:=(m_{1},m_{1}')$ is a dec-min
element of $\odotZ{B}$ and hence $m_{1}$ is indeed the restriction of a
dec-min element of $\odotZ{B}$ to $S_{1}$.

(B2)$\rightarrow $(A1):
By Theorems \ref{optkrit} and \ref{smallest},
we have $m_{1}(s)\geq \beta_{1} -1$ for each $s\in S_{1}(m)=S_{1}$, that is,
$\beta_{1} -1\leq m_{1}(s)\leq \beta_{1}$.  
By letting $r'$ denote the number of $\beta_{1}$-valued components of $m_{1}$, 
we obtain by \eqref{(r11)} that 
\[
r_{1}+(\beta_{1} -1)\vert S_{1}\vert = p_{1}(S_{1})
    = \widetilde m_{1}(S_{1})  = (\beta_{1} -1)\vert S_{1}\vert +r'
\] 
and hence $r'=r_{1}$.  
\finbox

\medskip 
Theorem~\ref{felbont} implies that, in order to characterize
the set of dec-min elements of $\odotZ{B}$, it suffices to characterize
the set of dec-min elements of $\odotZ{B_{1}'}$.

\begin{theorem}   \label{beta2} 
Let $\beta_{2}$ denote the smallest integer for
which $\odotZ{B_{1}'}$ has a $\beta_{2}$-covered element, that is, 
$\beta_{2}= \beta (B_{1}')$.  
Then 
\begin{equation} 
\beta_{2} = \max \{ \llceil {p_{1}'(X) \over \vert X\vert }\rrceil :  
                  \emptyset \not =X\subseteq S-S_{1}\},
\label{(betamin2)} 
\end{equation}
 where $p_{1}'(X)= p(X\cup S_{1}) - p(S_{1})$.
Furthermore, $\beta_{2}$ is the largest component in $S-S_{1}$ of every
dec-min element of $\odotZ{B}$, and $\beta_{2} < \beta_{1}$.  
\end{theorem}  

\Proof 
Formula \eqref{(betamin2)} follows by applying Theorem~\ref{betamin} 
to base-polyhedron $B_{1}'$ ($=B'(p_{1}')$) in place of $B$.
By Theorem~\ref{felbont}, the largest component in $S-S_{1}$ of any
dec-min element $m$ of $\odotZ{B}$ is $\beta_{2}$.  
By Theorem~\ref{smallest}, $S_{1}(m)=S_{1}$, and the definition of $S_{1}(m)$ shows
that $m(s)\leq \beta_{1} -1$ holds for every $s\in S-S_{1}$, from which
$\beta_{2}<\beta_{1}$ follows.  
\finbox

\subsection{The matroid $M_{1}$ on $S_{1}$} 
\label{M1}

It is known from the theory of base-polyhedra that the intersection of
an integral base-polyhedron with an integral box is a (possibly empty)
integral base-polyhedron.  
Moreover, if the box in question is small, 
then the intersection is actually a translated matroid base-polyhedron 
(meaning that the intersection arises from a matroid base-polyhedron by
translating it with an integral vector).  
This result is a consequence of the theorem that 
($*$) any integral base-polyhedron in the unit
$(0,1)$-cube is the convex hull of (incidence vectors of) the bases of a matroid.

Consider the special small integral box
$T_{1}\subseteq {\bf Z}\sp{S_{1}}$ defined by 
\[
T_{1}:=\{x:  \beta_{1} -1\leq x(s)\leq \beta_{1}\}
\]
and its intersection $B_{1}\sp{\bullet} := B_{1}\cap T_{1}$ 
with the base-polyhedron
$B_{1}$ investigated above.  
Therefore $B_{1}\sp{\bullet}$ is a translated matroid
base-polyhedron and Theorem~\ref{m1.ekv} implies the following.

\begin{corollary}   \label{trans.mat} 
The dec-min elements of $\odotZ{B_{1}}$ are exactly 
the integral elements of the translated matroid base-polyhedron $B_{1}\sp{\bullet}$.  
\finbox 
\end{corollary}

Our next goal is to reprove Corollary \ref{trans.mat} by concretely
describing the matroid in question 
and not relying on the background theorem ($*$) mentioned above.  
For a dec-min element $m_{1}$ of $\odotZ{B_{1}}$, let 
\[
   L_{1}(m_{1}):= \{s\in S_{1}: m_{1}(s)=\beta_{1}\}.
\] 
\noindent
We know from Theorem~\ref{m1.ekv} that 
$\vert L_{1}(m_{1})\vert =r_{1}$.  
Define a set-system ${\cal B}_{1}$ as follows:
\begin{equation} 
{\cal B}_{1}:=\{ L\subseteq S_{1}:  L=L_{1}(m_{1})
 \ \hbox{ for some dec-min element $m_{1}$ of $\odotZ{B_{1}}\}.$}\ 
\label{(B1def)} 
\end{equation}

We need the following characterization of ${\cal B}_{1}$.

\begin{proposition} \label{PRbaseM1}
 An $r_{1}$-element subset $L$ of $S_{1}$ is in ${\cal B}_{1}$
if and only if 
\begin{equation} 
\vert L\cap X\vert \geq p_{1}'(X):= p_{1}(X) - (\beta_{1} -1)\vert X\vert  
\ \hbox{ whenever } \  X\subseteq S_{1}.  
\label{(p01)} 
\end{equation} 
\end{proposition}

\Proof 
Suppose first that $L\in {\cal B}_{1}$, that is, 
there is a dec-min element $m_{1}$ of $\odotZ{B_{1}}$ for which $L=L_{1}(m_{1})$.  
Then
\[
(\beta_{1} -1)\vert X\vert + \vert X\cap L\vert = \widetilde m_{1}(X) \geq p_{1}(X),
\] 
for every subset $X\subseteq S_{1}$ from which \eqref{(p01)} follows.

To see the converse, let $L\subseteq S_{1}$ be an $r_{1}$-element set meeting \eqref{(p01)}.  
Let
\begin{equation} 
m_{1}(s):= \begin{cases} 
          \beta_{1} & \ \ \hbox{if}\ \ \ s\in L \cr
         \beta_{1} -1 & \ \ \hbox{if}\ \ \ s\in S-L.  
         \end{cases} 
\end{equation}
\noindent
Then obviously $L=L_{1}(m_{1})$.  Furthermore,
\[
\widetilde m_{1}(S_{1}) 
  = (\beta_{1} -1)\vert S_{1}\vert + \vert L\vert 
  = (\beta_{1} -1)\vert S_{1}\vert + r_{1} =p(S_{1})
\] 
and
\[
  \widetilde m_{1}(X) 
   = (\beta_{1} -1)\vert X\vert + \vert L\cap X\vert
  \geq p_{1}(X) \ \hbox{whenever}\ \ X\subset S_{1},
\] 
showing that $m_{1}\in B_{1}$.  
Since $m_{1}\in T_{1}$, we conclude that $m_{1}$ is a dec-min element of $\odotZ{B_{1}}$.  
\finbox

\begin{theorem}   \label{matroid1} 
The set-system \ ${\cal B}_{1}$ defined in \eqref{(B1def)} 
forms the set of bases of a matroid $M_{1}$ on ground-set $S_{1}$.  
\end{theorem}

\Proof 
The set-system ${\cal B}_{1}$ is clearly non-empty and all of its
members are of cardinality $r_{1}$.
It is widely known \cite{Edmonds70} 
that for an integral submodular function $b$ on a ground-set $S_{1}$ the set-system 
\[
\{ \  L\subseteq S_{1}:  \ \vert L\cap X\vert \leq b(X) \ \hbox{ whenever } \ 
 X\subset S_{1}, \ \vert L\vert = b(S_{1}) \ \},
\] 
if non-empty,
satisfies the matroid basis axioms.  This implies for the supermodular
function $p_{1}'$ that the set-system 
$\{L:  \vert L\cap X\vert \geq p_{1}'(X) \ \hbox{ whenever } \ 
     X\subset S_{1}, \ \vert L\vert = p_{1}'(S_{1}) \}$, 
if non-empty, forms the set of bases of a matroid.  
By applying this fact to the supermodular function $p_{1}'$ defined by 
$p_{1}'(X):= p_{1}(X) - (\beta_{1} -1)\vert X\vert $, 
one obtains that ${\cal B}_{1}$ is non-empty
and forms the set of bases of a matroid.  
\finbox

\medskip

With this matroid $M_{1}$, we can rewrite Corollary \ref{trans.mat} 
into a more explicit form,
which is convenient for our subsequent discussion.

\begin{corollary}   \label{M1B1} 
Let $\Delta_{1}:S_{1}\rightarrow {\bf Z} $ denote
the integral vector defined by $\Delta_{1}(s):=\beta_{1} -1$  for $s\in S_{1}$.  
A member $m_{1}$ of $\odotZ{B_{1}}$ is decreasingly minimal 
if and only if 
there is a basis $B_{1}$ of $M_{1}$ such that $m_{1}=\chi_{B_{1}} + \Delta_{1}$.
\finbox 
\end{corollary}

\subsection{Value-fixed elements of $S_{1}$}
\label{value-fixed}

We say that an element $s\in S$ is {\bf value-fixed} with respect to $\odotZ{B}$ 
if $m(s)$ is the same for every dec-min element $m$ of $\odotZ{B}$.  
In Section~\ref{SCoptdualset}, we will show a description of value-fixed
elements of $\odotZ{B}$.  In the present section, we consider the
value-fixed elements with respect to $B_{1}$, 
that is, $s\in S_{1}$ is value-fixed if 
$m_{1}(s)$ is the same for every dec-min element $m_{1}\in \odotZ{B_{1}}$.  
Recall that $m_{1}\in \odotZ{B_{1}}$ was shown to be dec-min precisely if 
$\beta_{1}-1\leq m_{1}(s)\leq \beta_{1}$ for each $s\in S_{1}$.

A {\bf loop} of a matroid is an element $s\in S_{1}$ not belonging to any basis.  
(Often the singleton $\{s\}$ is called a loop, that is,
$\{s\}$ is a one-element circuit).  
A {\bf co-loop} (or cut-element or isthmus) 
of a matroid is an element $s$ belonging to all bases.

\begin{proposition}  
$M_{1}$ has no loops.  
\end{proposition}  

\Proof
 By Proposition \ref{PR=S1}, for every $s\in S_{1}$ there is a
pre-dec-min element $m$ of $\odotZ{B}$ for which $m(s)=\beta_{1}$.  
Then $m_{1}:=m\vert S_{1}$ is a pre-dec-min element of $\odotZ{B_{1}}$ 
by Theorem~\ref{m1.ekv} from which $s_{1}$ belongs to a basis of $M_{1}$ 
by Corollary \ref{M1B1}.  
\finbox

\medskip

The proposition implies that:

\begin{proposition}   \label{fixbetai} 
If $s\in S_{1}$ is value-fixed (with respect to $B_{1}$), 
then $m_{1}(s)=\beta_{1}$ for every dec-min element $m_{1}$ of $\odotZ{B_{1}}$.  
\finbox 
\end{proposition}

By Corollary \ref{M1B1}, an element $s\in S_{1}$ is a co-loop of $M_{1}$
if and only if 
$m_{1}(s)=\beta_{1}$ holds for every dec-min element $m_{1}$ of $\odotZ{B_{1}}$.  
This and Theorem~\ref{felbont} imply the following.

\begin{theorem}   \label{co-loop} 
For an element $s\in S_{1}$, the following properties are pairwise equivalent. 
\smallskip

\noindent 
{\rm (A)} \ $s$ is a co-loop of $M_{1}$.  
\smallskip

\noindent 
{\rm (B)} \ $s$ is value-fixed.  
\smallskip

\noindent 
{\rm (C)} \ $m(s)=\beta_{1}$ holds for every dec-min element $m$ of $\odotZ{B}$.  
\finbox 
\end{theorem}

Our next goal is to characterize the set of value-fixed elements of $S_{1}$.  
Consider the family of subsets $S_{1}$ defined by 
\begin{equation} 
{\cal F}_{1}:= \{X\subseteq S_{1}:  \ \beta_{1} \vert X\vert = p_{1}(X)\}.
\label{(scriptF1)} 
\end{equation}
The empty set belongs to ${\cal F}_{1}$ and it is possible that 
${\cal F}_{1}$ has no other members.  
By standard submodularity arguments,
${\cal F}_{1}$ is closed under taking union and intersection.  
Let $F_{1}$ denote the unique largest member of ${\cal F}_{1}$.  
It is possible that $F_{1}=S_{1}$ in which case we call $S_{1}$ {\bf degenerate}.

\begin{theorem}   \label{F1} 
An element $s\in S_{1}$ is value-fixed 
if and only if 
$s\in F_{1}$.  
\end{theorem}  

\Proof 
Let $m_{1}$ be a dec-min member of $\odotZ{B_{1}}$.  Then
\[
 \beta_{1}\vert F_{1}\vert \geq \widetilde m_{1}(F_{1}) \geq p_{1}(F_{1}) =
\beta_{1}\vert F_{1}\vert 
\]
 and hence we must have $\beta_{1}=m_{1}(s)$ for
every $s\in F_{1}$, that is, the elements of $F_{1}$ are indeed value-fixed.

Conversely, let $s$ be value-fixed, that is, $m_{1}(s)=\beta_{1}$ for
each dec-min element $m_{1}$ of $\odotZ{B_{1}}$.  
Let $m_{1}$ be a dec-min member of $\odotZ{B_{1}}$.  
Let $Z$ denote the unique smallest set containing $s$
for which $\widetilde m_{1}(Z)= p_{1}(Z)$.  
(That is, $Z=T_{m_{1}}(s;p_{1}).)$
We claim that $m_{1}(t)=\beta_{1}$ for every element $t\in Z$.  
For if $m_{1}(t)=\beta_{1} -1$ for some $t$, then 
$m'_{1}:=m_{1}- \chi_{s}+\chi_{t}$ 
would also be a dec-min member of $\odotZ{B_{1}}$, 
contradicting the assumption that $s$ is value-fixed.
Therefore $ p_{1}(Z) = \widetilde m_{1}(Z) = \beta_{1}\vert Z\vert $ from
which the definition of $F_{1}$ implies 
that $Z\subseteq F_{1}$ and hence $s\in F_{1}$.  
\finbox

\section{Description of the set of all decreasingly minimal elements}
\label{canonical}

Let $B=B'(p)$ denote again an integral base-polyhedron defined by the
(integer-valued) supermodular function $p$.  
As in the previous section, 
$\odotZ{B}$ continues to denote the M-convex set consisting of
the integral vectors (points, elements) of $B$.  
Our present goal is
to provide a complete description of the set of all decreasingly minimal
(= egalitarian) elements of $\odotZ{B}$
by identifying a partition of the ground-set, 
to be named the canonical partition,
inherent in this problem.
As a consequence, we show that the set of
dec-min elements has a matroidal structure and this feature makes it
possible to solve the minimum cost dec-min problem.

\subsection{Canonical partition and canonical chain}

In Section~\ref{peak} we introduced the integer $\beta_{1}$ as the
minimum of the largest component of the elements of $\odotZ{B}$ 
as well as the peak-set $S_{1}$.
The peak-set $S_{1}$ induces a face of $B$, which is 
the direct sum of base-polyhedra
$B_{1}=B'(p_{1})$ and $B_{1}'=B'(p_{1}')$, 
where $p_{1}$ denotes the restriction
of $p$ to $S_{1}$ and $p_{1}'$ 
is obtained from $p$ by contracting $S_{1}$ 
(that is, $p_{1}'(X)=p(S_{1}\cup X)-p(S_{1})$).

A consequence of Theorem~\ref{felbont} is that, in order to
characterize the set of dec-min elements of $\odotZ{B}$, 
it suffices to characterize separately the dec-min elements 
of $\odotZ{B_{1}}$ and the dec-min elements of $\odotZ{B_{1}'}$.  
By Theorem~\ref{m1.ekv}, the dec-min elements of $\odotZ{B_{1}}$ are characterized 
as those elements of $\odotZ{B_{1}}$ which belong to the small box 
$T_{1}:=\{x\in {\bf R}\sp{S_{1}}:  \beta_{1} -1\leq x(s)\leq \beta_{1}$
 \ for \ $s\in S_{1}\}$.  
If the peak-set $S_{1}$ happens to be
the whole ground-set $S$, then the characterization of the set of
dec-min elements of $\odotZ{B}$ is complete.  
If $S_{1}\subset S$, then
our remaining task is to characterize the set of dec-min elements of
$\odotZ{B_{1}'}$.  This can be done by repeating iteratively the
separation procedure to the base-polyhedron 
$B_{1}' = B'(p_{1}') \subseteq {\bf R}\sp{S-S_{1}}$ 
described in Section~\ref{peak} for $B$.

In this iterative way, we are going to define a partition 
${\cal P}\sp{*}=\{S_{1},S_{2}, \dots , S_{q}\}$ of $S$ 
which determines a chain 
${\cal C}\sp{*}=\{C_{1},C_{2},\dots ,C_{q}\}$ where 
$C_{i}:=S_{1}\cup S_{2}\cup \cdots \cup S_{i}$ (in particular $C_{q}=S)$, 
and the supermodular function 
\[
p_{i}':=p/C_{i} \quad  \mbox{ on set } \  \overline{C_{i}}:=S-C_{i} 
\]
which defines the base-polyhedron 
$B_{i}'=B'(p_{i}')$ in ${\bf R}\sp{\overline{C_{i}}}$.  
Moreover, we define iteratively a decreasing sequence 
$\beta_{1}>\beta_{2}>\cdots >\beta_{q}$ 
of integers, a small box 
\begin{equation} 
 T_{i}:=\{x\in {\bf R}\sp{S_{i}}: \ \beta_{i}-1 \leq x(s)\leq \beta_{i} \ \hbox{for}\ s\in S_{i}\},
\label{(boxi)} 
\end{equation} 
and the supermodular function $p_{i}$ on $S_{i}$, where 
\begin{equation} 
p_{i}:=p_{i-1}'\vert S_{i}  \quad (=(p/C_{i-1})\vert S_{i}), 
\label{(pidef)} 
\end{equation} 
that is,
\[
 p_{i}(X)= p(X\cup C_{i-1}) - p(C_{i-1}) \ \hbox{for}\ \ X\subseteq S_{i}. 
\]
\noindent
Let $B_{i} := B'(p_{i})\subseteq {\bf R}\sp{S_{i}}$
be the base-polyhedron defined by $p_{i}$.

In the general step, suppose that the pairwise disjoint non-empty sets
$S_{1}, S_{2},\dots ,S_{j-1}$ 
have already been defined, along with the decreasing sequence 
$\beta_{1}>\beta_{2}>\cdots >\beta_{j-1}$ 
of integers.  
If $S=S_{1}\cup \cdots \cup S_{j-1}$, then by taking $q:=j-1$, 
the iterative procedure terminates.  
So suppose that this is not the case, that is, $C_{j-1}\subset S$.  
We assume that $p_{j-1}$
on $S_{j-1}$ has been defined as well as $p_{j-1}'$ on $\overline{C_{j-1}}$.

Let 
\begin{equation} 
\beta_{j}=\max \{ \llceil {p_{j-1}'(X) \over \vert X\vert }\rrceil :  
    \emptyset \not =X\subseteq \overline{C_{j-1}}\},
\label{(betaminj} 
\end{equation}
that is,
\begin{equation} 
\beta_{j}=\max \{ \llceil {p(X\cup C_{j-1}) - p(C_{j-1}) \over \vert X\vert }\rrceil :  
    \emptyset \not =X\subseteq \overline{C_{j-1}}\}.
\label{(betaj} 
\end{equation}
\noindent
Note that, by the iterative feature of these definitions,
Theorem~\ref{beta2} implies that 
\[
\beta_{j}<\beta_{j-1}.  
\]
\noindent
Furthermore, let $h_{j}$ be a set-function on $\overline{C_{j-1}}$ defined as follows:
\begin{equation} 
h_{j}(X):= p_{j-1}'(X) - (\beta_{j}-1)\vert X\vert 
 \ \ \hbox{for }\ \ X\subseteq \overline{C_{j-1}}, 
\label{(hj.def)}
\end{equation}
and let $S_{j}\subseteq \overline{C_{j-1}}$ be 
the peak-set of $\overline{C_{j-1}}$ assigned to $B_{j-1}':= B'(p_{j-1}')$,
that is, $S_{j}$ is the smallest subset of 
$\overline{C_{j-1}}$ maximizing $h_{j}$.  
Finally, let $p_{j}:=p_{j-1}'\vert S_{j}$ and let $p_{j}':= p_{j-1}'/S_{j}$.  
Observe by \eqref{(Z1Z2)} that $p_{j}'= p/C_{j}$.  
Therefore $p_{j}$ is a set-function on $S_{j}$ while $p_{j}'$ is defined on $\overline{C_{j}}$.

We shall refer to the partition ${\cal P}\sp{*}$ and the chain 
${\cal C}\sp{*}$ defined above as the {\bf canonical partition} and 
{\bf canonical chain} of $S$, respectively, assigned to $B$, 
while the sequence $\beta_{1}>\cdots >\beta_{q}$ 
will be called the {\bf essential value-sequence} of $\odotZ{B}$.

\begin{example} \rm \label{EXori4cano}
In the orientation problem in Example \ref{EXori4}, 
the essential value-sequence is given by $\beta_{1}=3$ and $\beta_{2}=2$ (with $q=2$)
and the canonical partition is  
${\cal P}\sp{*} = \{ S_{1}, S_{2} \}$
with 
$S_{1} = \{ a,b \}$ and $S_{2} = \{ c,d \}$.
The canonical partition for the M-convex set arising from such orientation problem
coincides with the ``density decomposition'' of Borradaile et al.~\cite{BMW}. 
\finbox
\end{example}

\begin{example} \rm \label{EXlongestVScano}
In the proof of Theorem \ref{equi.1},
we considered the longest chain consisting 
of non-empty $m$-tight and $m$-top sets, 
where $m$ is a dec-min element of the M-convex set $\odotZ{B}$.  
One may wonder  
whether this longest chain is the same as the canonical chain of $\odotZ{B}$.  
The following example, from the area of graph orientations, 
demonstrates that the answer is negative.

Let $G=(V,E)$ be an undirected graph where $V=\{a,b, c,d, x,y\}$ and
$E$ consists of four parallel edges between $a$ and $b$, four parallel edges
between $c$ and $d$, and the following four further edges:  $ax$, $cx$, $by$, $dy$.  
Consider the M-convex set $\odotZ{B}$ consisting of the in-degree vectors of 
all possible orientations of $G$.  
Graph $G$ has an orientation where the in-degree of every node is 2. 
This shows that the vector $m:=(2,2,2,2,2,2)$ is in $\odotZ{B}$,
and this uniform vector is obviously a dec-min element of $\odotZ{B}$.
For this $m$, the chain 
$\{ a,b \} \subset \{ a,b,c,d \} \subset \{ a,b,c,d,x \} \subset V$ 
is a longest chain consisting of four non-empty $m$-tight and $m$-top sets 
(and there are three other longest chains).
But the canonical chain of $\odotZ{B}$ consists of the single member $\{ V \}$.  
\finbox
\end{example}

Let $B\sp{\oplus}$ denote the face of $B$ 
defined by the canonical chain ${\cal C}\sp{*}$, 
that is, $B\sp{\oplus}$ is the direct sum of the $q$ base-polyhedra 
$B'(p_{i}) \ (i=1,\dots ,q)$.  
Finally, let $T\sp{*}$ be the direct sum of the small boxes 
$T_{i}$ $(i=1,\dots ,q)$, 
that is, $T\sp{*}$ is the integral box
defined by the essential value-sequence as follows:
\begin{equation} 
T\sp{*}:= \{x\in {\bf R}\sp{S}:  \ \beta_{i}-1\leq x(s)\leq \beta_{i} 
  \ \ \hbox{whenever}\ s\in S_{i} \ (i=1,\dots ,q)\}, 
\label{(box)} 
\end{equation}
and let 
\[
B\sp{\bullet}:= B\sp{\oplus}\cap T\sp{*}.
\]
This set $B\sp{\bullet}$ is an integral base-polyhedron,
since the intersection of an integral base-polyhedron with an integral box 
is always an integral base-polyhedron.
Furthermore, $B\sp{\bullet}$ is the direct sum of the $q$ base-polyhedra 
$B_{i}\cap T_{i}$ ($i=1,\dots ,q$), where $B_{i}=B'(p_{i})$, 
implying that a vector $m$ is in $\odotZ{B\sp{\bullet}}$ 
if and only if each $m_{i}$ is in $\odotZ{B_{i}}\cap T_{i}$, where $m_{i}=m\vert S_{i}$.

\begin{theorem}  \label{main.1} 
Let $B=B'(p)$ be an integral base-polyhedron on ground-set $S$.  
The set of decreasingly minimal elements of $\odotZ{B}$ 
is (the {\rm M}-convex set) $\odotZ{B\sp{\bullet}}$.
Equivalently, an element $m\in \odotZ{B}$ is decreasingly minimal 
if and only if 
its restriction $m_{i}:=m\vert S_{i}$ to $S_{i}$ belongs to
$B_{i}\cap T_{i}$ for each $i=1,\dots ,q$, 
where $\{S_{1},\dots ,S_{q}\}$ is the canonical partition of $S$ belonging to $B$, 
\ $T_{i}$ is the small box defined in \eqref{(boxi)}, 
and $B_{i}$ is the base-polyhedron
$B'(p_{i})$ belonging to the supermodular set-function $p_{i}$ defined in \eqref{(pidef)}.  
\end{theorem}

\Proof 
We use induction on $q$.  
Suppose first that $q=1$, that is,  $S_{1}=S$ and $B_{1}=B$.  
If $m$ is a dec-min element of $B$, 
then the equivalence of Properties (A1) and (A3) in Theorem~\ref{m1.ekv}
implies that $m$ is in $\odotZ{B\sp{\bullet}}$.  
If, conversely, $m\in \odotZ{B\sp{\bullet}}$, 
then $m$ is near-uniform and, by the equivalence of
Properties (A1) and (A3) in Theorem~\ref{m1.ekv} again, $m$ is dec-min.

Suppose now that $q\geq 2$ and consider the base-polyhedron 
$B_{1}' = B'(p_{1}')$ 
appearing in Theorem~\ref{felbont}.  
The iterative definition of the canonical partition ${\cal P}\sp{*}$ 
implies that the canonical partition of $S-S_{1}$ assigned to 
$B_{1}'$ is $\{S_{2},\dots ,S_{q}\}$ 
and the essential value-sequence belonging to $B_{1}'$ is 
$\beta_{2}>\beta_{3}>\cdots >\beta_{q}$.  
Also, the canonical chain 
${\cal C}':= \{C_{2}',\dots ,C_{q}'\}$
of $B_{1}'$ consists of the sets 
$C_{i}'=S_{2}\cup \cdots \cup S_{i} = C_{i}-S_{1}$ \ $(i=2,\dots ,q)$.

By applying the inductive hypothesis to $B_{1}'$, we obtain that an
integral element $m_{1}'$ of $B_{1}'$ is dec-min 
if and only if
$m_{1}'$ is in the face of $B_{1}'$ defined by chain ${\cal C}'$ and $m_{1}'$ 
belongs to the box 
$T':= \{x\in {\bf R}\sp{S-S_{1}}:  \ \beta_{i}-1\leq x(s)\leq
   \beta_{i} \ \ \hbox{whenever}\ s\in S_{i} \ (i=2,\dots ,q)\}$.  
By applying Theorem~\ref{felbont}, we are done in this case as well.
\finbox 
\medskip

\begin{corollary} \label{COchardecmin}
Let $B=B'(p)$ be an integral base-polyhedron on ground-set $S$.  
Let $\{C_{1},\dots ,C_{q}\}$ be the canonical chain,
$\{S_{1},\dots ,S_{q}\}$ the canonical partition of $S$,
and $\beta_{1}> \beta_{2}> \dots > \beta_{q}$ 
the essential value-sequence belonging to $\odotZ{B}$.  
Then an element $m\in \odotZ{B}$ is decreasingly minimal 
if and only if 
each $C_{i}$ is $m$-tight 
(that is, $\widetilde m(C_{i})=p(C_{i})$) 
and $\beta_{i} -1 \leq m(s) \leq \beta_{i}$ holds for each 
$s\in S_{i}$ \ $(i=1,\dots ,q)$.  
\finbox 
\end{corollary}

\subsection{Obtaining the canonical chain and value-sequence from a dec-min element}

The main goal of this section is to show that 
the canonical chain and value-sequence can be rather easily obtained 
from an arbitrary dec-min element of $\odotZ{B}$.
This approach will be crucial in developing a polynomial algorithm in \cite{FM19partB}
for computing the essential value-sequence along with  the canonical chain and partition.

Let $m$ be an element of $\odotZ{B}$.  
We called a set $X\subseteq S$ $m$-tight if $\widetilde m(X)=p(X)$.  
Recall from Section~\ref{egal2base} that, 
for a subset $Z\subseteq S$, $T_{m}(Z)=T_{m}(Z;p)$
denoted the unique smallest $m$-tight set including $Z$, that is,
$T_{m}(Z)$ is the intersection of all the $m$-tight sets including $Z$.
Obviously,
\begin{equation} 
T_{m}(Z)= \cup (T_{m}(z) :  z\in Z). 
 \label{(Tm)} 
\end{equation}

Let $m$ be an arbitrary dec-min element of $\odotZ{B}$.  
We proved that $m$ is in the face $B\sp{\oplus}$ of $B$ defined 
by the canonical chain 
${\cal C}\sp{*} =\{C_{1},\dots ,C_{q}\}$
belonging to $B$.  
Therefore each $C_{i}$ is $m$-tight with respect to $p$.  
Furthermore $m_{i}:=m\vert S_{i}$ 
belongs to the box $T_{i}$ defined in \eqref{(boxi)}.  
This implies that 
$m(s)\geq \beta_{i}-1$ for every $s\in C_{i}$ and 
$m(s')\leq \beta_{i+1}$ for every $s'\in \overline{C_{i}}$. 
 (The last inequality holds indeed
since $s'\in \overline{C_{i}}$ implies that $s'\in S_{j}$ for some $j\geq i+1$
from which $m(s')\leq \beta_{j}\leq \beta_{i+1}$.)  
Since $\beta_{i+1}\leq \beta_{i}-1$, we obtain that each $C_{i}$ is an $m$-top set.

Since $m_{i}$ is near-uniform on $S_{i}$ with values $\beta_{i}$ and
possibly $\beta_{i}-1$, we obtain
\[
 \beta_{i}= \llceil {\widetilde m_{i} (S_{i}) \over \vert S_{i}\vert }\rrceil
= \llceil {p_{i} (S_{i}) \over \vert S_{i}\vert }\rrceil 
= \llceil {p (C_{i})- p(C_{i-1}) \over \vert S_{i}\vert }\rrceil .
\]
\noindent
Let 
$L_{i}:=\{s\in S-C_{i-1}:  m(s) =\beta_{i}\}$
and let $r_{i}:=\vert L_{i}\vert $. 
Then $p_{i}(S_{i}) = \widetilde m_{i}(S_{i}) = (\beta_{i}-1)\vert
S_{i}\vert + r_{i}$ and hence
\begin{equation} 
r_{i}= p(C_{i}) - p(C_{i-1}) - (\beta_{i}-1)\vert S_{i}\vert .
\label{(ri)} 
\end{equation}

The content of the next lemma is that, once $C_{i-1}$ is given, the
next member $C_{i}$ of the canonical chain (and hence $S_{i}$, as well)
can be expressed with the help of $m$.  
Recall that
$T_{m}(L_{i})=T_{m}(L_{i};p)$ denoted the smallest $m$-tight set including $L_{i}$.

\begin{lemma}   \label{Sim} 
$C_{i}= C_{i-1} \cup T_{m}(L_{i};p)$.  
\end{lemma}

\Proof 
Recall the definition of function $h_{i}$ given in \eqref{(hj.def)}.  
We have 
\begin{equation} 
h_{i}(S_{i}) = r_{i} 
\label{(hiSi)} 
\end{equation}
since 
$h_{i}(S_{i}) 
  =p_{i-1}'(S_{i}) -(\beta_{i}-1)\vert S_{i}\vert 
  = p(S_{i}\cup C_{i-1}) - p(C_{i-1}) - (\beta_{i}-1)\vert S_{i}\vert 
  =\widetilde m(C_{i}) - \widetilde m(C_{i-1}) - (\beta_{i}-1)\vert S_{i}\vert
  = \widetilde m(S_{i}) - (\beta_{i}-1)\vert S_{i}\vert = r_{i}$.

Since $L_{i}\subseteq C_{i}$ and each of $C_{i-1}$, $C_{i}$, and $T_{m}(L_{i})$
are $m$-tight, we have $C_{i-1} \cup T_{m}(L_{i};p)\subseteq C_{i}$.  
For $X':= T_{m}(L_{i}) \cap \overline{C_{i-1}}$ we have
\begin{eqnarray*} h_{i}(X') 
&=& p(C_{i-1}\cup T_{m}(L_{i})) -p(C_{i-1}) - (\beta_{i}-1)\vert X_{i}'\vert 
\\ &=& 
\widetilde m(C_{i-1}\cup T_{m}(L_{i}))
- \widetilde m(C_{i-1}) - (\beta_{i}-1)\vert X_{i}'\vert 
\\ &=&
\widetilde m(X') - (\beta_{i}-1)\vert X_{i}'\vert = \vert L_{i}\vert =r_{i} = h_{i}(S_{i}), 
\end{eqnarray*}
that is, $X'$ is also a maximizer of $h_{i}(X)$.  
Since $S_{i}$ was the smallest maximizer of $h_{i}$, 
we conclude that $C_{i-1} \cup T_{m}(L_{i};p) \supseteq C_{i}$.  
\finbox

\medskip 
The lemma implies that both the essential value-sequence
$\beta_{1}>\cdots >\beta_{q}$ 
and the canonical chain ${\cal C}\sp{*}$
belonging to $\odotZ{B}$ can be directly obtained from $m$.

\begin{corollary}   \label{Cim} 
Let $m$ be an arbitrary dec-min element of $\odotZ{B}$.  
The essential value-sequence and the canonical chain belonging to $\odotZ{B}$ 
can be described as follows.
Value $\beta_{1}$ is the largest $m$-value and $C_{1}$ is the smallest
$m$-tight set containing all $\beta_{1}$-valued elements.  
Moreover, for $i=2,\dots ,q,$ \ $\beta_{i}$ is the largest value of 
$m\vert \overline{C_{i-1}}$ and $C_{i}$ is the smallest $m$-tight set 
(with respect to $p$) 
containing each element of $m$-value at least $\beta_{i}$.
\finbox 
\end{corollary}

A detailed algorithm based on this corollary will be described in
\cite{FM19partB}.
Note that a dec-min element $m$ of $\odotZ{B}$ 
may have more than $q$ distinct values.  
For example, if $q=1$ and $L_{1}\subset C_{1}=S$, then
$m$ has two distinct values, namely $\beta_{1}$ on the elements of
$L_{1}$ and $\beta_{1} -1$ on the elements of $S-L_{1}$,
while its essential value-sequence consists of the single member $\beta_{1}$.

\medskip

\paragraph{A direct proof} 
Corollary \ref{Cim} implies that
the chain of subsets and value-sequence assigned to a dec-min element $m$ of $\odotZ{B}$ 
in the corollary do not depend on the choice of $m$.
Here we describe an alternative, direct proof of this consequence.

\begin{theorem} \label{THnew54} 
Let $m$ be an arbitrary dec-min element of $\odotZ{B}$.  
Let $\beta_{1}$ denote the largest value of $m$ and 
let $C_{1}$ denote the smallest $m$-tight set (with respect to $p$)
containing all $\beta_{1}$-valued elements.  
Moreover, for $i=2,3,\dots,q$, let $\beta_{i}$ denote 
the largest value of $m\vert \overline{C_{i-1}}$
and let $C_{i}$ denote the smallest $m$-tight set 
containing each element of $m$-value at least $\beta_{i}$.  
Then the chain $C_{1}\subset C_{2}\subset \cdots \subset C_{q}$ 
and the sequence 
$\beta_{1}>\beta_{2}>\cdots >\beta_{q}$ do not depend on the choice of $m$.
\end{theorem}

\Proof 
Let $z$ be dec-min element of $\odotZ{B}$.  
We use induction on the number of elements $t$ of $S$ for which $m(t)>z(t)$.  
If no such an element $t$ exists, then $m=z$ and there is nothing to prove.  
So assume that $z\not =m$.

Let $L_{i}:=\{t\in S_{i}:  m(t)=\beta _{i}\}$.  
As $m$ is dec-min, the definition of $C_{i}$ implies that 
$m(s)=\beta_{i}-1$ holds for every element 
$s \in S_{i}-L_{i}$.
Let $t\in L_{i}$ and let $s\in T_{m}(t)-L_{i}$.
Then $m':= m + \chi_{s} - \chi_{t}$ is also a dec-min element of $\odotZ{B}$, 
and we say that $m'$ is obtained from $m$ by an elementary step.  
Observe that $T_{m}(t) = T_{m'}(s)$ and hence the chain and 
the value-sequence assigned to $m'$ is the same as those assigned to $m$.

Let $i$ denote the smallest subscript for which $m\vert S_{i}$ and $z\vert S_{i}$ differ.  
Since $z$ is dec-min, $z(s)\leq \beta _{i}$ holds for every $s\in S_{i}$.  
Let $L'_{i}:=\{t\in S_{i}:  z(t)=\beta _{i}\}$.  
Then
$z(v)\leq \beta _{i}-1$ for every $v\in S_{i}-L'_{i}$, 
and $\vert L'_{i}\vert \leq \vert L_{i}\vert $ as $z$ is dec-min.  
Therefore
\[
\widetilde z(S_{i})
 \ \leq \ 
\beta _{i}\vert L'_{i}\vert + (\beta _{i}-1)(\vert S_{i}-L'_{i}\vert ) 
\ = \ 
(\beta _{i}-1)\vert S_{i}\vert + \vert L'_{i}\vert 
\ \leq \
(\beta _{i}-1)\vert S_{i}\vert + \vert L_{i}\vert .
\]
On the other hand,
\begin{align*}
\widetilde z(S_{i}) & = \widetilde z(C_{i}) -\widetilde z(C_{i-1}) =
\widetilde z(C_{i}) - \widetilde m(C_{i-1}) 
\\
 & \geq p(C_{i}) -
\widetilde m(C_{i-1}) = \widetilde m(C_{i}) - \widetilde m(C_{i-1})=
\widetilde m(S_{i})= (\beta _{i}-1)\vert S_{i}\vert + \vert L_{i}\vert .
\end{align*}
\noindent
Therefore we have equality throughout, in particular, 
$\widetilde z(C_{i})=p(C_{i})$, 
$\vert L'_{i}\vert = \vert L_{i}\vert $, and $z(v)=\beta_{i}-1$ for every $v\in S_{i}-L'_{i}$.

Let $t\in L_{i}$ be an element for which $m(t)>z(t)$.  
Then $m(t)=\beta_{i}$ and $z(t)=\beta _{i}-1$.  
It follows that $T_{m}(t)$ contains an element $s$ for which $z(s) > m(s)$, 
implying that $m(s)=\beta _{i}-1$ and $z(s)=\beta _{i}$.  
Now $m(t)>m'(t) = z(t)$ holds for the dec-min element 
$m':= m + \chi_{s} - \chi_{t}$ 
obtained from $m$ by an elementary step, and therefore we are done by induction.  
\finbox

\subsection{Matroidal description of the set of dec-min elements}
\label{SCmatrdescdmset}

In Section~\ref{M1}, we introduced a matroid $M_{1}$ on $S_{1}$ and proved in 
Corollary \ref{trans.mat}
that the dec-min elements of $\odotZ{B_{1}}$
are exactly the integral elements of the translated base-polyhedron of $M_{1}$, 
where the translation means the addition of the constant vector
$(\beta_{1} -1,\dots ,\beta_{1} -1)$ of dimension $\vert S_{1}\vert $. 
The same notions and results can be applied to each subscript $i=2,\dots,q$.
Furthermore, by formulating Lemma~\ref{matroid1} 
for subscript $i$ in place of 1, we obtain the following.

\begin{proposition} \label{matroidi} 
The set-system \
${\cal B}_{i}:=\{ L\subseteq S_{i}:  L=L_{i}(m_{i})$ 
{\rm for some dec-min element} $m_{i}$ {\rm of} $\odotZ{B_{i}}\}$ \
forms the set of bases of a matroid $M_{i}$ on ground-set $S_{i}$.  
An $r_{i}$-element subset $L$ of $S_{i}$ is a basis of $M_{i}$ if and only if 
\begin{equation} 
\vert L\cap X\vert  \geq    p_{i}'(X)  
  := p_{i}(X)  - (\beta_{i}-1)\vert X\vert 
\label{(p0)} 
\end{equation} 
holds for every $X\subseteq S_{i}$.
\finbox 
\end{proposition}

It follows that a vector $m_{i}$ on $S_{i}$ is a dec-min element of $\odotZ{B_{i}}$ 
if and only if 
$\beta _{i}-1\leq m_{i}(s)\leq \beta _{i}$ 
for each $s\in S_{i}$ and the set 
$L_{i}:=\{s\in S_{i}:  m_{i}(s)=\beta _{i}\}$ is a basis of $M_{i}$.  
Let $M\sp{*}$ denote the direct sum of matroids $M_{1},\dots ,M_{q}$ and
let $\Delta \sp{*}\in {\bf Z}\sp{S}$ denote the translation vector
defined by 
\[
\Delta \sp{*}(s):= \beta_{i}-1 \ \ \hbox{whenever}\ \ s\in S_{i}, \ i=1,\dots ,q.
\]

By integrating these results, we obtain the following characterization.

\begin{theorem}   \label{matroid-eltolt} 
Let $B$ be an integral base-polyhedron.  
An element $m$ of (the M-convex set) $\odotZ{B}$ is decreasingly minimal 
if and only if 
$m$ can be obtained in the form $m=\chi_{L}+ \Delta \sp{*}$
where $L$ is a basis of the matroid $M\sp{*}$.  
The base-polyhedron $B\sp{\bullet}$ 
arises from the base-polyhedron of $M\sp{*}$ 
by adding the translation vector $\Delta\sp{*}$.  
Concisely, the set of dec-min elements of $\odotZ{B}$ is a
matroidal M-convex set.  
\finbox 
\end{theorem}

\paragraph{Cheapest dec-min element} 
An important algorithmic consequence of Theorems \ref{main.1} and \ref{matroid-eltolt} is 
that they help solve the 
{\bf cheapest dec-min element} problem, 
which is as follows.  
Let $c:S\rightarrow {\bf R} $ be a cost function and
consider the problem of computing a dec-min element $m$ of an M-convex
set $\odotZ{B}$ for which $cm$ is as small as possible.

By Theorem~\ref{matroid-eltolt} the set 
$\odotZ{B\sp{\bullet}}$ 
of dec-min elements of
$\odotZ{B}$ can be obtained from a matroid $M\sp{*}$ 
by translation.
Namely, there is a vector $\Delta \sp{*} \in {\bf Z}\sp{S}$ 
such that $m$ is in 
$\odotZ{B\sp{\bullet}}$ 
if and only if 
there is a basis $L$ of $M\sp{*}$
for which  $m= \chi_{L} + \Delta \sp{*}$.  
Note that the matroid $M\sp{*}$
arises as the direct sum of matroids $M_{i}$ defined on the members
$S_{i}$ of the canonical partition.  
$M_{1}$ is described in Proposition \ref{PRbaseM1}  
and the other matroids $M_{i}$ may be determined analogously in an iterative way.  
To realize this algorithmically, 
we must have a strongly polynomial algorithm 
to compute the canonical partition 
as well as the essential value-sequence.  
Such an algorithm will be described in \cite{FM19partB}.

Therefore, in order to find a minimum $c$-cost dec-min element of $\odotZ{B}$, 
it suffices to find a minimum $c$-cost basis of $M\sp{*}$.
Note that, in applying 
the greedy algorithm to the matroids $M_{i}$
in question, 
we need a rank oracle, which can be realized with the help
of a submodular function minimization oracle 
by relying on the definition of bases in \eqref{(p01)}.

Recall that for integral bounds $f \leq g$,
the intersection $B_{1}$ of a base-polyhedron $B$
and the box $T(f,g)$, if non-empty, 
is itself a base-polyhedron.
Therefore the algorithm above can be applied to the 
M-convex set $\odotZ{B_{1}}$, that is, 
we can compute a cheapest dec-min element of the intersection
$\odotZ{B_{1}} = \odotZ{B} \cap T(f,g)$.

\section{Integral square-sum and difference-sum minimization}
\label{SCsqdiffsum}

For a vector  $z \in {\bf Z}\sp{S}$, 
we can conceive several natural functions
to measure the uniformity of its component values $z(s)$ for $s \in S$.
Here are two examples:
\begin{align}
&\mbox{\bf square-sum}: \quad \ \  
W(z):=\sum [z(s)\sp{2}:s\in S],
\label{squaresumdef}
\\
&\mbox{\bf difference-sum}: \  
\Delta(z):=\sum [ | z(s) - z(t) |:s \not= t, \  s,t \in S].
\label{diffsumdef}
\end{align}
\noindent
For vectors  $z_{1}$ and $z_{2}$ with
$\widetilde z_{1}(S)=\widetilde z_{2}(S)$,
$z_{1}$ may be felt more uniform than $z_{2}$ if $W(z_{1}) < W(z_{2})$,
and $z_{1}$ may also be felt more uniform 
if $\Delta(z_{1}) < \Delta(z_{2})$.
The first goal of this section is to show,
by establishing a fairly general theorem,
that a dec-min element of an M-convex set $\odotZ{B}$
is simultaneously a minimizer of these two functions.
The second goal of this section is to derive a
min-max formula for the minimum integral square-sum of an element of an M-convex set 
$\odotZ{B}$, along with characterizations of (integral) square-sum minimizers
and dual optimal solutions.

\subsection{Symmetric convex minimization} 
\label{SCnonsepmin}

Let $S$ be a non-empty ground-set of $n$ elements: 
$S = \{ 1,2,\ldots,n \}$.
We say that a function $\Phi: \ZZ\sp{S} \to \RR$ is {\bf symmetric} if
\begin{equation} \label{symmtry}
 \Phi(z(1), z(2), \ldots, z(n)) = 
 \Phi(z(\sigma(1)), z(\sigma(2)), \ldots, z(\sigma(n))) 
\end{equation}
for all permutations $\sigma$ of $(1,2,\ldots,n)$.
We call a function $\Phi: \ZZ\sp{S} \to \RR$ {\bf convex} if
\begin{equation} \label{nonstrcnvI}
 \lambda \Phi(x) + (1-\lambda) \Phi(y) \geq  \Phi( \lambda x + (1-\lambda) y)
\end{equation}
whenever
$x, y \in \ZZ\sp{S}$, \ $0 < \lambda < 1$, and  
$\lambda x + (1-\lambda) y$ is an integral vector;
and {\bf strictly convex} if
\begin{eqnarray}  \label{strcnvI}
 \lambda \Phi(x) + (1-\lambda) \Phi(y) >  \Phi( \lambda x + (1-\lambda) y)
\end{eqnarray}
whenever
$x, y \in \ZZ\sp{S}$, \ $0 < \lambda < 1$, 
and  $\lambda x + (1-\lambda)  y$ is an integral vector.

In the special case 
of a function in one variable,
it can easily be shown that the convexity of 
$\varphi: \ZZ \to \RR$
is equivalent to the weaker requirement that the inequality
\begin{equation}  \label{univarconvfndef}
  2\varphi (k) \leq \varphi (k-1) + \varphi (k+1)
\end{equation}
holds for every integer $k$.
Such function $\varphi$ is often called a (univariate) discrete convex function.
It is strictly convex in the sense of \eqref{strcnvI} if and only if 
$2\varphi (k) < \varphi (k-1) + \varphi (k+1)$ 
holds for every integer $k$.
For example, $\varphi(k)=k\sp{2}$ is strictly convex 
while $\varphi (k)=\vert k \vert$ is convex but not strictly.  
Given a function $\varphi $ in one variable
satisfying \eqref{univarconvfndef},
define $\Phi$ by 
\begin{equation}  \label{symsepconvfndef}
\Phi(z) := \sum [\varphi (z(s)):  s\in S]
\end{equation}
for $z \in {\bf Z}\sp{S}$.  
Such a function $\Phi$ is called a {\bf symmetric separable convex function};
note that $\Phi$ is indeed convex in the sense of \eqref{nonstrcnvI}.
When $\varphi $ is strictly convex, $\Phi$ is also strictly convex.

\begin{example} \rm \label{RMsquaresum}
The square-sum $W(z)$ in \eqref{squaresumdef}
is a symmetric convex function
which is separable and strictly convex.
\finbox
\end{example}

\begin{example} \rm \label{EXdiffsum}
The difference-sum $\Delta(z)$ in \eqref{diffsumdef}
is a symmetric convex function
which is neither separable nor strictly convex.
More generally, for a nonnegative integer $K$, 
the function defined by
\[
\Delta_K(z) := \sum [ ( |z(s) - z(t)|  - K)\sp{+}:s \not= t, \  s,t \in S]
\]
is a symmetric convex function,
where $(x )\sp{+} = \max \{ x, 0 \}$.
\finbox
\end{example}

The following statements show a close relationship 
between decreasing minimality and 
the minimization of symmetric convex $\Phi$ over an M-convex set $\odotZ{B}$.

\begin{proposition} \label{decminIsPhiMin=nonsep} 
Let $B$ be an integral base-polyhedron
and $\Phi$ a symmetric convex function.  
Then each dec-min element of $\odotZ{B}$ 
is a minimizer of $\Phi$ over $\odotZ{B}$.  
\end{proposition} 

\Proof 
Since the dec-min elements of $\odotZ{B}$
are value-equivalent and $\Phi$ is symmetric,
the $\Phi$-value of each dec-min element is the same value $\mu$.
We claim that $\Phi(m)\geq \mu$ for each $m \in \odotZ{B}$.  
Suppose indirectly that there is an element $m$ of 
$\odotZ{B}$ for which $\Phi(m) < \mu$.
Then $m$ is not dec-min in $\odotZ{B}$
and Property (A) in Theorem~\ref{equi.1} implies that 
there is a 1-tightening step for $m$
resulting in decreasingly smaller member of $\odotZ{B}$,
that is,
there exist $s, t \in S$ such that
$m(t) \geq m(s)  + 2$
and 
$m' := m + \chi_{s} - \chi_{t} \in \odotZ{B}$.

Let $\alpha = m(t) - m(s)$, where $\alpha \geq 2$,
and define
$z = m + \alpha (\chi_{s} - \chi_{t})$.
Since $z$ is obtained from $m$ by 
interchanging the components at $s$ and $t$,
the symmetry of $\Phi$ formulated in  (\ref{symmtry}) implies that $\Phi(m)=\Phi(z)$.  
Note that the vector $z$ may not be a member of $\odotZ{B}$.
For $\lambda = 1 - 1/\alpha$ we have
\begin{equation}  \label{decminnonsep=pr1}
 \lambda m + (1-\lambda)  z = 
\left( 1 - \frac{1}{\alpha} \right) m 
+ \frac{1}{\alpha} \left( m + \alpha (\chi_{s} - \chi_{t}) \right)  
=  m + \chi_{s} - \chi_{t} = m'
\in \odotZ{B} \ (\subseteq \ZZ^{S} ) ,
\end{equation}
from which 
$\lambda \Phi(m) + (1-\lambda) \Phi(z) \geq  \Phi(m')$
by convexity (\ref{nonstrcnvI}).
Since $\Phi(m) = \Phi(z)$,
this implies
$\Phi(m) \geq  \Phi(m')$.
After a finite number of
such 1-tightening steps, we arrive at a dec-min element
$m_{0}$ of $\odotZ{B}$,
for which $\mu = \Phi(m_{0})\leq \Phi(m) < \mu$, a contradiction.
\finbox

\medskip 

Note that if $\Phi $ is convex but not strictly convex,
then $\Phi$ may have minimizers that are not dec-min elements.  
This is exemplified by the identically zero function $\Phi$ for which
every member of $\odotZ{B}$ is a minimizer.  
However, for strictly convex functions we have the following characterization.

\begin{theorem}  \label{decmin=Phi=nonsep} 
Given an integral base-polyhedron $B$ and
a symmetric strictly convex function $\Phi$, 
an element $m$ of $\odotZ{B}$ is a minimizer of $\Phi$ 
if and only if 
$m$ is a dec-min element of $\odotZ{B}$.  
\end{theorem}  

\Proof
 If $m$ is a dec-min element, then $m$ is a $\Phi$-minimizer by
Proposition \ref{decminIsPhiMin=nonsep}.  
To see the converse, let $m$ be a
$\Phi$-minimizer of $\odotZ{B}$.  If, indirectly, $m$ is not a dec-min
element, then Property (A) in Theorem~\ref{equi.1} implies that 
there is a 1-tightening step for $m$,
that is,
there exist $s, t \in S$ such that
$m(t) \geq m(s)  + 2$
and 
$m' := m + \chi_{s} - \chi_{t} \in \odotZ{B}$.
For $\alpha = m(t) - m(s)$,
$\lambda = 1 - 1/\alpha$, and
$z = m + \alpha (\chi_{s} - \chi_{t})$,
we have
\eqref{decminnonsep=pr1},
from which we obtain
$\lambda \Phi(m) + (1-\lambda) \Phi(z) >  \Phi(m')$
by strict convexity (\ref{strcnvI}).
Since $z$ is obtained from $m$ by interchanging the components at $s$ and $t$, 
the symmetry of $\Phi$ formulated in  (\ref{symmtry}) implies that $\Phi(m)=\Phi(z)$.  
But then $\Phi(m)>\Phi(m')$ would follow, 
in contradiction to the assumption that $m$ is a minimizer of $\Phi$.
\finbox

\medskip

We obtain the following  as corollaries of this theorem.

\begin{corollary}  \label{COdecmin=separable} 
Let $B$ be an integral base-polyhedron
and $\Phi$ a symmetric separable convex function.  
Then each dec-min element of $\odotZ{B}$ 
is a minimizer of $\Phi$ over $\odotZ{B}$,
and the converse is also true 
if, in addition,  $\Phi$ is strictly  convex.
\finbox
\end{corollary}  

\begin{corollary}  \label{COdecmin=squaresum} 
For an M-convex set $\odotZ{B}$,
an element $m$ of $\odotZ{B}$ is a square-sum minimizer 
if and only if 
$m$ is a dec-min element of $\odotZ{B}$.  
\finbox
\end{corollary}

An immediate consequence of Corollary \ref{COdecmin=separable}
is that a square-sum
minimizer of $\odotZ{B}$ minimizes an arbitrary symmetric 
separable convex function $\Phi$.
Note, however, that this consequence immediately follows from a much earlier result 
of Groenevelt \cite{Groenevelt91} below, which deals with
the minimization of a (not-necessarily symmetric) separable convex function.

\begin{theorem}[Groenevelt \cite{Groenevelt91}; cf.~{\cite[Theorem 8.1]{Fujishigebook}}]
 \label{THgroenevelt}
Let $B$ be an integral base-polyhedron,
$\odotZ{B}$ be the set of its integral elements,
and $\Phi(z) = \sum [\varphi_{s}(z(s)):  s\in S]$ \ for $z\in \ZZ\sp{S}$, 
where $\varphi_{s}: \ZZ \to \RR \cup \{ +\infty \}$
 is a discrete convex function for each $s\in S$.
An element $m$ of $\odotZ{B}$ is a minimizer of $\Phi(z)$
if and only if
$\varphi_{s}(m(s)+1) + \varphi_{t}(m(t)-1) \geq \varphi_{s}(m(s)) + \varphi_{t}(m(t))$
whenever 
$m+\chi_{s}-\chi_{t} \in \odotZ{B}$.
\finbox
\end{theorem}

A dec-min element is also characterized as a difference-sum minimizer.

\begin{theorem}  \label{THdecmin=diffsum} 
For an M-convex set $\odotZ{B}$,
an element $m$ of $\odotZ{B}$ is a difference-sum minimizer 
if and only if $m$ is a dec-min element of $\odotZ{B}$.  
\end{theorem}

\Proof
By Proposition \ref{decminIsPhiMin=nonsep}
every dec-min element is a difference-sum minimizer.
To show the converse, suppose indirectly that there is 
difference-sum minimizer $m$
that is not dec-min in $\odotZ{B}$.
Property (A) in Theorem~\ref{equi.1} implies that 
there is a 1-tightening step for $m$,
that is,
there exist $s, t \in S$ such that
$m(t) \geq m(s)  + 2$
and 
$m' := m + \chi_{s} - \chi_{t} \in \odotZ{B}$.
Here we observe that
$|m'(s) - m'(t)| =  |m(s) - m(t)|  -2$
and
\[
(|m'(v) - m'(s)| + |m'(v) - m'(t)|) - (|m(v) - m(s)| + |m(v) - m(t)|)
=
\begin{cases}
 -2 &  \hbox{if}\ \ m(s) < m(v) < m(t)
\\
  \phantom{-} 0 &  \hbox{otherwise}.
\end{cases} 
\]
This shows $\Delta(m') \leq \Delta(m) -2$, a contradiction.
\finbox

\medskip

\begin{remark} \rm \label{RMorderingsqsum}
Corollary \ref{COdecmin=squaresum} says 
that an element $m$ of an M-convex set
$\odotZ{B}$ is dec-min precisely if $m$ is a square-sum minimizer.
One may feel that it would
have been a more natural approach to derive this equivalence by
showing that $x \leq_{\rm dec} y$ holds precisely if $W(x)\leq W(y)$.
Perhaps surprisingly, however, this equivalence fails to hold, that
is, the square-sum is not order-preserving with respect to the
quasi-order $\leq_{\rm dec}$.  
To see this, consider the following
four vectors in increasing order:
\[
 m_{1}=(2,3,3,1) <_{\rm dec} m_{2}=(3,3,3,0) <_{\rm dec} m_{3}=(2,2,4,1)
 <_{\rm dec} m_{4}=(3,2,4,0).
\]
\noindent
Their square-sums admit a different order:
\[
W(m_{1})= 23, \quad W(m_{2})=27, \quad  W(m_{3}) = 25, \quad W(m_{4}) = 29.
\]
The four vectors $m_{i}$ $(i=1,2,3,4)$ form an M-convex set.
Among these four elements, $m_{1}$
is the unique dec-min element and the unique square-sum minimizer but
the decreasing-order and the square-sum order of the other three
elements are different.
We remark that if $\varphi$
in \eqref{symsepconvfndef}
is not only strictly convex but \lq rapidly\rq \ increasing 
as well, then $x <_{\rm dec} y$
can be proved to be equivalent to $\Phi(x) < \Phi (y)$.  
This intuitive notion of rapid increase is formalized in 
\cite{FM19partII}.
\finbox
\end{remark}

\begin{remark} \rm \label{RMintersqsum}
For the intersection of two M-convex sets,
dec-min elements and square-sum minimizers may not coincide.  
Here is an example.
Let $\odotZ{B_{1}}=\{ (3,3,3,0 ), \ (2,2,4,1), \ (2,3,3,1), \ (3,2,4,0) \}$
and 
$\odotZ{B_{2}}=\{ (3,3,3,0 ), \ (2,2,4,1), \ (3,2,3,1), \ (2,3,4,0) \}$,
which are both M-convex.
In their intersection
$\odotZ{B_{1}} \cap \odotZ{B_{2}} = \{ (3,3,3,0), (2,2,4,1) \}$,
the vector $(3,3,3,0)$ is the unique dec-min element
while $(2,2,4,1)$ is the unique square-sum minimizer.
This demonstrates
that the two notions of optima may differ 
for the intersection of two M-convex sets.
\finbox
\end{remark}

\begin{remark} \rm \label{RMklargesymconv}
In Section~\ref{SCklargestsum}, we considered the minimum $k$-largest-sum problem
that aimed at finding an element of $\odotZ{B}$ 
for which  the sum of the $k$ largest components is as small as possible.
For each $k$, the sum of the $k$ largest components is a symmetric convex function, 
and hence Theorem \ref{klargest} is a special case 
of Proposition \ref{decminIsPhiMin=nonsep}.
It is noted, however, that
Theorem \ref{decmin=Phi=nonsep} is not applicable to the $k$-largest-sum problem,
as this function is not strictly convex.
Nevertheless, a dec-min element can be characterized in terms of 
the $k$-largest-sum
if we simultaneously consider the functions for all $k$
(Theorem \ref{THmajordecmin}).
\finbox
\end{remark}

\begin{remark} \rm \label{RMtwosepconv}
For $a,b,c \geq 0$, the function defined by
\[
\Phi(z) = a \sum_{s \in S} |z(s)| 
  + b \sum_{s \neq t} |z(s) - z(t)| + c \sum_{s \neq t} |z(s)+z(t)| 
\] 
is a symmetric convex function.
More generally, a function of the form 
\[
\Phi(z) = \sum_{s \in S} \varphi_{1}(z(s)) 
 + \sum_{s \neq t} \varphi_{2} (|z(s) - z(t)|)
 + \sum_{s \neq t} \varphi_{3}(z(s)+z(t))  ,
\] 
where $\varphi_{1}, \varphi_{2}, \varphi_{3}: \ZZ \to \RR$
are
 (discrete) convex functions
(as defined in \eqref{univarconvfndef}),
is a symmetric convex function which is not separable.
Such a function is an example of the so-called 
2-separable convex functions.
By Theorem~\ref{decmin=Phi=nonsep}, 
a dec-min element of $\odotZ{B}$ is a minimizer of function $\Phi$ over $\odotZ{B}$.
The minimization of 2-separable convex functions 
is investigated in depth by Hochbaum and others \cite{AHO04,Hoc02,Hoc07}
using network flow techniques.
\finbox
\end{remark}

\begin{remark} \rm \label{RMnoduality}
A min-max formula  can be derived for the square-sum
(see Section~\ref{SCsquare.minmax})
and, more generally, for separable convex functions
from the Fenchel-type duality theorem in DCA
\cite{Murota98a,Murota03}.
However, we cannot use the Fenchel-type duality theorem 
to obtain a min-max formula for non-separable symmetric convex functions,
since non-separable symmetric convex functions are not necessarily {\rm M}-convex.
\finbox
\end{remark}

\subsection{Min-max theorem for integral square-sum}
\label{SCsquare.minmax}

Recall the notation
$W(z)=\sum [z(s)\sp{2}:s\in S]$
for the square-sum of $z \in {\bf Z}\sp{S}$.  
Given a polyhedron $B$, 
we say that an element $m\in \odotZ{B}$ is a {\bf square-sum minimizer} 
(over $\odotZ{B}$) or that $m$ is an {\bf integral square-sum minimizer} 
of $B$ 
if $W(m)\leq W(z)$ holds for each $z\in \odotZ{B}$.  
The main goal of this section is to derive a
min-max formula for the minimum integral square-sum of an element of an M-convex set 
$\odotZ{B}$, along with a characterization of (integral) square-sum minimizers.

A set-function $p$ on $S$ can be considered as a function defined on $(0,1)$-vectors.  
It is known that $p$ can be extended in a natural
way to every vector $\pi$ in ${\bf R}\sp{S}$, as follows.  
For the sake of this definition, we may assume that the elements of $S$ are
indexed in a decreasing order of the components of $\pi$, that is,
$\pi(s_{1})\geq \cdots \geq \pi(s_{n})$ (where the order of the
components of $\pi$ with the same value is arbitrary).
For $j=1,\dots ,n$,  let $I_{j}:=\{s_{1},\dots ,s_{j}\}$ and let
\begin{equation}
\hat p(\pi):= p(I_{n})\pi(s_{n}) 
 + \sum_{j=1}\sp {n-1} p(I_{j})[\pi(s_{j})-\pi(s_{j+1})].
\end{equation}
\noindent
Obviously, $p(Z)=\hat p(\chi_Z)$.  
The function $\hat p$ is called 
\cite[Section 14.5.1]{Frank-book}
the {\bf linear extension} of $p$,
where $\hat p$ is a piecewise-linear function in $\pi$.

\begin{remark} \rm \label{RMlovext}
The linear extension was first
considered by Edmonds \cite{Edmonds70} who proved for a polymatroid
$P=P(b)$ defined by a monotone, non-decreasing submodular function $b$ that 
$\max \{\pi x :  x\in \odotZ{P}\} = \hat b(\pi )$ when $\pi$ 
is non-negative.  
The same approach shows for a base-polyhedron $B=B'(p)$
defined by a supermodular function $p$ that 
$\min \{\pi x :  x\in \odotZ{B}\} = \hat p(\pi)$.  
Another basic result is due to Lov\'asz \cite{Lovasz83} 
who proved that $p$ is supermodular 
if and only if 
$\hat p$ is concave.
We do not, however, explicitly need these results, 
and only remark that in the literature the linear extension is often called
Lov\'asz extension.
\finbox
\end{remark}

\medskip

Our approach is as follows.  First, we consider an arbitrary
set-function $p$ on $S$ (supermodular or not) along with the
polyhedron
\[
 B=B'(p):=\{x:  x\in {\bf R}\sp{S}, \ \widetilde x(Z)\geq p(Z)
 \ \  \mbox{for every $Z\subset S$ and $\widetilde x(S)=p(S)$} \}, 
\]
and develop an easily checkable lower bound for the minimum
square-sum over the integral elements of $B$.  If this lower bound is
attained by an element $m$ of $\odotZ{B}$, then $m$ is certainly a
square-sum minimizer independently of any particular property of $p$.
For general $p$, the lower bound (not surprisingly) is not 
always attainable.  
We shall prove, however, that it is attainable when $p$ is supermodular.  
That is, we will have a min-max theorem for
the minimum square-sum over an M-convex set $\odotZ{B}$, or in other
words, we will have an easily checkable certificate for an element $m$
of $\odotZ{B}$ to be a minimizer of the square-sum.

We shall need the following two claims.
For any real number $\alpha \in {\bf R}$, 
let $\lfloor \alpha \rfloor $ denote the largest integer not larger than $\alpha$, 
and $\lceil \alpha \rceil$ the smallest integer not smaller than $\alpha$.

\begin{claim}   \label{hcestim}
 For $m,\pi \in {\bf Z}\sp{S}$, one has 
\begin{equation} 
\sum_{s\in S} \llfloor {\pi(s) \over 2}\rrfloor \llceil {\pi(s) \over 2}\rrceil 
    \geq \sum_{s\in S} m(s)[\pi(s)-m(s)].  
\label{(hcestim.1)}
\end{equation}
\noindent
Moreover, equality holds 
if and only if 
\begin{equation} 
 m(s)\in \bigg\{ \llfloor {\pi(s) \over 2}\rrfloor , 
                 \llceil {\pi(s) \over 2}\rrceil \bigg\}
   \quad \hbox{\rm for every} \ s\in S.
\label{(hcestim.2)} 
\end{equation} 
\end{claim}

\Proof 
The claim follows by observing that 
$\lfloor {a / 2}\rfloor \lceil {a  / 2}\rceil \geq b(a-b)$ 
holds for any pair of integers
$a$ and $b$, where equality holds precisely if 
$b \in \big\{\lfloor {a / 2}\rfloor , \lceil {a / 2}\rceil \big\}$.  
\finbox

\medskip

Let $p$ be an arbitrary set-function on $S$ with $p(\emptyset )=0$ and
consider an integral element $m$ of the polyhedron $B=B'(p)$.  
Recall that a non-empty subset $X\subseteq S$ was called 
a strict $\pi$-top set if $\pi(u)>\pi(v)$ held whenever $u\in X$ and $v\in S-X$.  
In what follows, for an $m\in \odotZ{B}$, $m$-tightness of a subset
$Z\subseteq S$ means $\widetilde m(Z)=p(Z)$.

\begin{claim}   \label{pimpi} 
For $m\in \odotZ{B}$ and $\pi \in {\bf Z}\sp{S}$, one has 
\begin{equation} 
\hat p(\pi) \leq \sum_{s\in S}m(s)\pi(s).
\label{(pimpi)} 
\end{equation} 
\noindent
Moreover, equality holds 
if and only if 
each (of the at most $n$) strict $\pi$-top set is $m$-tight.  
\end{claim}

\Proof 
Suppose that the elements of $S$ are indexed in such a way that
$\pi(s_{1})\geq \pi(s_{2})\geq \cdots \geq \pi(s_{n})$.  
For $j=1,\dots ,n$,   let $I_{j}:=\{s_{1},\dots ,s_{j}\}$.   
Then
\begin{eqnarray*} 
\hat p(\pi) &=& p(I_{n})\pi(s_{n}) 
    + \sum_{j=1}\sp{n-1} p(I_{j})[\pi(s_{j})-\pi(s_{j+1})] 
\\ &\leq & 
\widetilde m(I_{n})\pi(s_{n}) 
   + \sum_{j=1}\sp{n-1} \widetilde m(I_{j})[\pi(s_{j})-\pi(s_{j+1})]
\\ &=& 
\sum_{1\leq i\leq j\leq n} m(s_{i})\pi(s_{j}) 
   - \sum_{1\leq i\leq j\leq n-1} m(s_{i})\pi(s_{j+1}) 
\\ &=&
 \sum_{1\leq i\leq j\leq n} m(s_{i})\pi(s_{j})  
  - \sum_{1\leq i< j'\leq n} m(s_{i})\pi(s_{j'}) 
\\ &=& 
\sum_{j=1}\sp {n} m(s_{j})\pi(s_{j}), 
\end{eqnarray*}
from which \eqref{(pimpi)} follows.  
Furthermore, we have equality in \eqref{(pimpi)} 
precisely if $\widetilde m(I_{j})=p(I_{j})$
holds whenever $\pi(s_{j})-\pi(s_{j+1})>0$.  
But this latter condition is equivalent to requiring 
that each strict $\pi$-top set is $m$-tight.  
\finbox 
\medskip

\begin{proposition}   \label{lowerbound}
Let $p$ be an arbitrary set-function on $S$ with $p(\emptyset )=0$ 
and let $m$ be an integral element of the polyhedron $B=B'(p)$.
Then
\begin{equation} 
 \sum_{s\in S} m(s)\sp{2}
 \geq \hat p(\pi) - \sum_{s\in S} \llfloor {\pi(s) \over 2}\rrfloor 
 \llceil {\pi(s) \over 2}\rrceil 
\label{(minleqmax)} 
\end{equation}
whenever $\pi \in {\bf Z}\sp{S}$ is an integral vector.
Furthermore, equality holds for $m$ and $\pi$ 
if and only if 
the following {\bf optimality criteria} hold:
\begin{eqnarray} 
{\rm (O1)} & & \hbox{\rm \eqref{(hcestim.2)} holds:} \ \
  m(s)\in 
  \bigg\{ \llfloor {\pi(s) \over 2}\rrfloor , \llceil {\pi(s) \over 2}\rrceil \bigg\} 
 \quad \hbox{\rm for every }\  s\in S , 
\label{(optcrit.1)} 
\\ 
{\rm (O2)} & & \hbox{\rm each strict $\pi$-top-set is $m$-tight with respect to $p$}.  
\label{(optcrit.2)}
\end{eqnarray}
\end{proposition}  

\Proof 
Let $\pi \in {\bf Z}\sp{S}$.  By the two preceding claims,
\begin{equation} 
\sum_{s\in S} m(s)\sp{2} 
   = \sum_{s\in S} m(s)\pi(s) 
   - \sum_{s\in S} m(s) [\pi(s)-m(s)] 
   \geq \hat p(\pi) 
   - \sum_{s\in S} \llfloor {\pi(s) \over 2}\rrfloor \llceil {\pi(s) \over 2}\rrceil ,
\label{(westim)} 
\end{equation} 
from which \eqref{(minleqmax)} follows.
The claims also immediately imply that we have equality in
\eqref{(minleqmax)} precisely if the optimality criteria (O1) and (O2) hold.  
\finbox
\medskip

The min-max formula in the next theorem concerning min square-sum over
the integral elements of an integral base-polyhedron can be derived from 
the more general Fenchel-type duality theorem in DCA 
(see \cite{Murota98a} and also Theorem 8.21, page 222, in the book \cite{Murota03}),
or from  a recent framework 
\cite{FM20boxTDI} 
of separable discrete convex function minimization over 
the integer points in an integral box-TDI polyhedron.
However, our proof relies only on the relatively simple
characterization of dec-min elements described in Theorem~\ref{equi.1}.  
In particular, we need no results of Sections \ref{peak} and \ref{canonical}.

\begin{theorem}   \label{minmax}
Let $B=B'(p)$ be a base-polyhedron defined by
an integer-valued supermodular function $p$.  
Then
\begin{equation} 
 \min \{ \sum_{s\in S} m(s)\sp{2}:  m\in \odotZ{B} \} 
 = \max \{\hat p(\pi) - \sum_{s\in S} \llfloor {\pi(s) \over 2}\rrfloor 
 \llceil {\pi(s) \over 2}\rrceil :  \pi \in {\bf Z}\sp{S} \}.  
\label{(minmax)} 
\end{equation}
\end{theorem}

\Proof 
By Proposition \ref{lowerbound}, 
$\min\geq \max$ holds in \eqref{(minmax)} and hence all what 
we have to prove is that there is
an element $m\in \odotZ{B}$ and an integral vector 
$\pi \in {\bf Z}\sp{S}$ meeting the two optimality criteria formulated 
in Proposition \ref{lowerbound}.
Let $m$ be an arbitrary dec-min element of $\odotZ{B}$.  
By Property
(B) of Theorem~\ref{equi.1}, there is a chain $(\emptyset \subset ) \
C_{1}\subset C_{2}\subset \cdots \subset C_\ell= S$ of $m$-tight and
$m$-top sets for which the restrictions of $m$ onto the difference sets 
$S_{i}:=C_{i}-C_{i-1}$ ($i=1,\dots ,\ell$) are near-uniform in $S_{i}$
(where $C_{0}:=\emptyset $).  
Note that $\{S_{1},\dots ,S_\ell\}$ is a partition of $S$.

For $i=1,\dots ,\ell$, let $\beta_{i}(m):= \max \{m(s):  s\in S_{i}\}$.
Define $\pi_{m}:S\rightarrow {\bf Z} $ by
\[
 \pi_{m}(s):= 2\beta_{i}(m)-1  \ \hbox{ if } \ s\in S_{i} \ \ (i=1,\dots ,\ell)  . 
\]
We have 
\[
\lfloor {\pi_{m}(s)/ 2}\rfloor 
  = \beta_{i}(m)-1 \leq m(s)
 \leq \beta_{i}(m)= \lceil {\pi_{m}(s)/2} \rceil 
\]
for every $s\in S_{i}$, and
hence Optimality criterion (O1) holds for $m$ and $\pi_{m}$.

We claim that each strict $\pi_{m}$-top set $Z$ is a member of chain $\cal C$.  
Indeed, as $\pi_{m}$ is uniform in each $S_{j}$, 
if $Z$ contains an element of $S_{j}$, then $Z$ includes the whole $S_{j}$.
Furthermore, since each member of $\cal C$ is an $m$-top set, we have
$\beta_{1}(m)\geq \beta_{2}(m)\geq \cdots \geq \beta_\ell(m)$, 
and hence if
$Z$ includes $S_{j}$, then it includes each $S_{i}$ with $i<j$.  
Therefore
every strict $\pi_{m}$-top set is indeed a member of the chain,
implying Optimality criterion (O2).  
\finbox

\medskip

It should be noted that the optimal dual solution $\pi_{m}$ obtained in
the proof of the theorem is actually an {\bf odd} vector in the sense
that each of its component is an odd integer.

\begin{corollary}   \label{minmaxb} 
There is an odd dual optimizer $\pi$ in
the min-max formula \eqref{(minmax)}, that is, the min-max formula in
Theorem {\rm \ref{minmax}} can be re-written as follows:
\begin{equation} 
\min \{ \sum_{s\in S} m(s)\sp{2}:  m\in \odotZ{B} \}
 = \max \{\hat p(\pi) - \sum_{s\in S} { {\pi(s)\sp{2} - 1} \over 4} : 
  \ \pi \in {\bf Z}\sp{S}, \ \pi \ \hbox{\rm is odd} \ \}.  
\label{(minmax.odd)} 
\end{equation}
\end{corollary}

We emphasize that for the proof of Theorem~\ref{minmax} and Corollary
\ref{minmaxb} we relied only on Theorem~\ref{equi.1} and did not need
the characterization of the set of dec-min elements of $\odotZ{B}$ 
given in Section~\ref{canonical}.

In the proof of Theorem~\ref{minmax}, we chose an arbitrary dec-min
element $m$ of $\odotZ{B}$ and an arbitrary chain of $m$-tight and
$m$-top sets 
such that $m$ is near-uniform on each difference set.  
In Section~\ref{canonical}, we proved that there is a single canonical
chain ${\cal C}\sp{*}$ which meets these properties for every dec-min
element of $\odotZ{B}$.  
Therefore the dual optimal $\pi\sp{*}$ assigned to ${\cal C}\sp{*}$ 
is also independent of $m$.  
Namely, consider the canonical $S$-partition $\{S_{1},\dots ,S_{q}\}$ and
the essential value-sequence $\beta_{1}>\cdots >\beta_{q}$.  
Define $\pi\sp{*}$ by
\begin{equation} 
\pi\sp{*}(s):= 2\beta_{i}-1 \ \hbox{ if } \ s\in S_{i} \ (i=1,\dots ,q).  
\label{(pi*def)} 
\end{equation}
\noindent
As we pointed out in the proof of Theorem~\ref{minmax}, this $\pi\sp{*}$ 
is also a dual optimum in \eqref{(minmax)}.  
We shall prove in the next section that 
$\pi\sp{*}$ is actually the unique smallest dual
optimum in \eqref{(minmax)}.

\subsection{The set of optimal duals to integral square-sum minimization}
\label{SCoptdualset}

We proved earlier that an element $m\in \odotZ{B}$ is a square-sum
minimizer precisely if it is a dec-min element.  
This and Theorem~\ref{main.1} imply that the square-sum minimizers 
of $\odotZ{B}$ are the integral members of a base-polyhedron $B\sp{\bullet}$
obtained by intersecting a particular face of $B$ with a special small box. 
This means that the integral square-sum minimizers form an M-convex set.

Our next goal is to reveal the structure of the set $\Pi$ of 
the dual optima in Theorem~\ref{minmax}
and we provide a description of $\Pi$ 
as the integral solution set of feasible potentials in a box.  
This shows another connection to DCA,
which is discussed after the proof of Theorem~\ref{minmax}.

Recall that the optimality criteria for a dec-min element $m$ of
$\odotZ{B}$ and for an integral vector $\pi$ were 
given by (O1) and (O2) in \eqref{(optcrit.1)}--\eqref{(optcrit.2)}. 
These immediately imply the following.

\begin{proposition}  
For an integral vector $\pi$, the following are equivalent.  
\smallskip

\noindent 
{\rm (A)} \ 
$\pi$ is a dual optimum (that is, $\pi$ belongs to $\Pi$).  
\smallskip

\noindent 
{\rm (B)} \ 
There is a dec-min element $m$ of $\odotZ{B}$
such that $m$ and $\pi$ meet the optimality criteria.
\smallskip

\noindent
{\rm (C)} \ 
For every dec-min $m$ of $\odotZ{B}$, 
$m$ and $\pi$ meet the optimality criteria.  
\finbox
\end{proposition}

Consider the canonical $S$-partition $\{S_{1},\dots ,S_{q}\}$, 
the essential value-sequence 
$\beta_{1}>\beta_{2}>\cdots >\beta_{q}$,
and the matroids $M_{i}$ on $S_{i}$ $(i=1,\dots ,q)$.  
We can use the notions and 
the results of Section~\ref{value-fixed} 
formulated for $M_{1}$ to each $M_{i}$ \ $(i=1,\dots ,q)$.  
To follow the pattern of ${\cal F}_{1}$ introduced in \eqref{(scriptF1)}, let 
\begin{equation}
 {\cal F}_{i}:= \{X\subseteq S_{i}:  \ \beta_{i}\vert X\vert = p_{i}(X)\} ,
\end{equation} 
where $p_{i}$ was defined by $p_{i}(X)=p(C_{i-1}\cup X) - p(C_{i-1})$ 
for $X\subseteq S_{i}$.  
Since $\beta_{i}\vert X\vert \geq p_{i}(X)$ for every $X\subseteq S_{i}$
and $p_{i}$ is supermodular, 
${\cal F}_{i}$ is closed under taking intersection and union.  
Let $F_{i}$ denote the unique largest member of ${\cal F}_{i}$, 
that is, $F_{i}$ is the union of the members of ${\cal F}_{i}$.  
Both $F_{i}=\emptyset $ and $F_{i}=S_{i}$ are possible.
As a generalization of Theorems \ref{co-loop} and \ref{F1} we obtain the following.

\begin{theorem}  \label{co-loopi} 
For an element $s\in S_{i}$ \ $(i=1,\dots ,q)$, 
the following properties are pairwise equivalent.%
\smallskip

\noindent 
{\rm (A)} \ $s$ is value-fixed.  
\smallskip

\noindent 
{\rm (B)} \ $m(s)=\beta_{i}$ holds for every dec-min element $m$ of $\odotZ{B}$.  
\smallskip

\noindent
 {\rm (C)} \ $s\in F_{i}$.  
\smallskip

\noindent
 {\rm (D)} \ $s$ is a co-loop of $M_{i}$.  
\finbox 
\end{theorem}

Define a digraph $D_{i}=(F_{i}, A_{i})$ on node-set $F_{i}$ in which $st$ is
an arc if $s,t\in F_{i}$ and there is no $t\overline{s}$-set in ${\cal F}_{i}$.
This implies that no arc of $D_{i}$ enters any member of ${\cal F}_{i}$.

\begin{theorem}  \label{opt.dual.set} 
An integral vector $\pi \in {\bf Z}\sp{S}$ 
is an optimal dual solution to the integral minimum square-sum problem 
(that is,  $\pi \in \Pi$) 
if and only if 
the following three conditions hold for each $i=1,\dots ,q:$
\begin{eqnarray} 
&& \pi(s)=2\beta_{i}-1  \quad \hbox{\rm  for every}\ \  s\in S_{i}-F_{i}, 
\label{(Si-Fi)} 
\\&& 
2\beta_{i}-1\leq \pi(s) \leq 2\beta_{i}+1 \quad \hbox{\rm for every}\ \ s\in F_{i}, 
\label{(Fi)} 
\\&& 
\pi(s)-\pi(t) \geq 0 \quad \hbox{\rm whenever \ $s,t\in F_{i}$ \  and \  $st\in A_{i}$}.  
\label{(pispit)} 
\end{eqnarray} 
\end{theorem}

\Proof
\begin{claim}  
\label{O1'} 
Optimality criterion {\rm (O1)} in \eqref{(optcrit.1)}
is equivalent to
\begin{equation} 
{\rm (O1')} \ \ \ \ \ 2m(s)-1 \leq \pi(s) \leq 2m(s)+1 \ \ 
\hbox{\rm for \ $s\in S$}. 
\label{(optcrit.1')} 
\end{equation} 
\end{claim} 

\Proof 
When $\pi(s)$ is even, we have the following equivalences:
\begin{eqnarray*} 
m(s)\in \bigg\{ \llfloor {\pi(s) \over 2}\rrfloor , 
           \llceil {\pi(s) \over 2}\rrceil \bigg\}
 & \Leftrightarrow & \pi(s)=2m(s) 
\\ &
\Leftrightarrow & 2m(s)-1 \leq \pi(s) \leq 2m(s)+1.  
\end{eqnarray*}
\noindent
When $\pi(s)$ is odd, we have the following equivalences:
\begin{eqnarray*} 
m(s)\in \bigg\{ \llfloor {\pi(s) \over 2}\rrfloor , 
                \llceil {\pi(s) \over 2}\rrceil \bigg\} 
& \Leftrightarrow & \pi(s)-1 \leq 2m(s) \leq \pi(s)+1 \ 
\\ & \Leftrightarrow & 2m(s)-1 \leq \pi(s) \leq 2m(s)+1.  
\end{eqnarray*}
\vspace{-2.8\baselineskip} \\
\finbox
\vspace{1.2\baselineskip}

Suppose first that $\pi \in {\bf Z}\sp{S}$ is an optimal dual solution.
Then the optimality criteria (O1$'$) and (O2) formulated in
\eqref{(optcrit.1')} and \eqref{(optcrit.2)} 
hold for every dec-min element $m$ of $\odotZ{B}$.

Let $s$ be an element of $S_{i}-F_{i}$.  
Since $s$ is not value-fixed, 
there are dec-min elements $m$ and $m'$ of $\odotZ{B}$ for which
$m(s)=\beta_{i}-1$ and $m'(s)=\beta_{i}$.  
By applying \eqref{(optcrit.1')} to $m$ and to $m'$, we obtain that
\[
 2\beta_{i}-1 = 2m'(s)-1\leq \pi(s)\leq 2m(s)+1 
   = 2(\beta_{i}-1) + 1 = 2\beta_{i}-1 ,
\]
from which $\pi(s)=2\beta_{i}-1$ follows, and hence \eqref{(Si-Fi)} holds indeed.

Let $s$ be an element of $F_{i}$.  
As $s$ is value-fixed, $m(s)=\beta_{i}$ holds 
for any dec-min element $m$ of $\odotZ{B}$.  
We obtain from \eqref{(optcrit.1')} that
\[
 2\beta_{i}-1 = 2m(s)-1\leq \pi(s)\leq 2m(s)+1 = 2\beta_{i} + 1 
\]
and hence \eqref{(Fi)} holds.

To derive \eqref{(pispit)}, suppose indirectly that $st$ is an arc in $A_{i}$ 
for which 
$\pi(t)>\pi(s) \ (\geq 2\beta_{i}-1)$.  
Let $Z:=\{v\in S: \pi(v)\geq \pi(t)\}$.  
Then $Z$ is a strict $\pi$-top set, and moreover, we have
$C_{i-1}\subseteq Z\subseteq (C_{i-1}\cup F_{i})-s$,
where the latter inclusion follows from \eqref{(Si-Fi)} and $\pi(t) > 2\beta_{i}-1$.  
By Optimality criterion (O2), $Z$ is $m$-tight with respect to $p$.  
Let $X:=Z\cap S_{i}$.  
Then $X\subseteq F_{i}$ and hence
\[
p(Z) = \widetilde m(Z) = \widetilde m(C_{i-1}) + \widetilde m(X) 
     = p(C_{i-1}) + \beta_{i}\vert X\vert , 
\]
from which 
\[
  \beta_{i}\vert X\vert = p(Z) - p(C_{i-1}) =p_{i}(X),
\]
that is, $X$ is in ${\cal F}_{i}$, 
contradicting
the definition of $A_{i}$ which requires that $st$ enters no member of ${\cal F}_{i}$.

Suppose now that $\pi$ meets the three properties 
\eqref{(Si-Fi)}, \eqref{(Fi)}, and \eqref{(pispit)}. 
Let $m\in \odotZ{B}$ be an arbitrary dec-min element.  
Consider an element $s$ of $S_{i}$.  
If $s\in F_{i}$, that is, if $s$ is value-fixed, then $m(s)=\beta_{i}$.  
By \eqref{(Fi)}, we have $2m(s)-1\leq \pi(s)\leq 2m(s)+1$, 
that is, Optimality criterion (O1$'$) holds.  
If $s\in S_{i}-F_{i}$, then $\pi(s)=2\beta_{i}-1$ by \eqref{(Si-Fi)}, 
from which 
\[
 \llfloor {\pi(s) \over 2}\rrfloor = {\pi(s)-1 \over 2} 
 \ = \  \beta_{i}-1
 \ \leq \  m(s) 
 \ \leq \  \beta_{i} 
 \ = \  {\pi(s) +1 \over 2} = \llceil {\pi(s) \over 2}\rrceil ,
\]
showing that Optimality criterion (O1$'$) holds.

To prove optimality criterion (O2), let $Z$ be a strict $\pi$-top set
and let $\mu := \min\{ \pi(v):  v\in Z\}$.  
 Let $i$ denote the largest subscript 
for which $X:=Z\cap S_{i}\not =\emptyset $. 
Then $\mu \leq 2\beta_{i}+1 \leq 2\beta_{i-1}-1\leq \pi(u)$ 
holds for every $u\in C_{i-1}$,  
from which $C_{i-1}\subseteq Z$ as $Z$ is a strict $\pi$-top set.

If $\mu =2\beta_{i}-1$, then $S_{i}\subseteq Z$ as $Z$ 
is a strict $\pi$-top set, from which 
$Z=C_{i}$, implying that $Z$ is an $m$-tight set in this case.  
Therefore we suppose 
$\mu \geq 2\beta_{i}$, from which $X\subseteq F_{i}$ follows.  
Now $X\in {\cal F}_{i}$, for otherwise there is an arc $st\in A_{i}$ $(s,t\in F_{i})$ 
entering $X$, and then $\pi(t)\leq \pi(s)$ holds by Property \eqref{(pispit)}; 
 this contradicts the assumption that $Z$ is a strict $\pi$-top set.  
By $X\in {\cal F}_{i}$ we have $\beta_{i}\vert X\vert =p_{i}(X)$,
whereas $m(s)=\beta_{i}$ for each $s \in X$ by $X \subseteq F_{i}$.  Hence
\begin{eqnarray*} 
\widetilde m(Z) &=& \widetilde m(X)+ \widetilde m(C_{i-1}) 
 = \beta_{i}\vert X\vert + p(C_{i-1}) 
\\ &=& p_{i}(X) + p(C_{i-1}) 
= p(X\cup C_{i-1}) - p(C_{i-1}) + p(C_{i-1}) = p(Z),
\end{eqnarray*}
that is, $Z$ is indeed $m$-tight.  
\finbox \finboxHere

\medskip \medskip

We now relate Theorem~\ref{minmax} to a concept from discrete convex analysis,
where two kinds of discrete convexity
play major roles as mutually `conjugate' notions of discrete convexity 
\cite{Murota98a,Murota03}.
One of them is  M-convexity and the other is called L-convexity.
One of the equivalent definitions says
 that a set $L$ of integer vectors is an
{\bf L-convex set}
if it is the set of integer-valued feasible potentials in the network flow problem.
Formally,
$L=\{\pi \in {\bf Z}\sp{S}:  \pi(v)-\pi(u) \leq g(uv) \ (u,v \in S) \}$,
where $g$ is an integer-valued function on the ordered pairs of elements of $S$.
A set of integer vectors is called an 
{\bf L$\sp{\natural}$-convex set}
(pronounce L-natural convex set) 
if it is the intersection of an L-convex set with an integral box.

In \eqref{(pi*def)}, we defined a special dual optimal solution 
$\pi\sp{*}$ by $\pi\sp{*}(s) =2\beta_{i}-1$ 
whenever $s\in S_{i}$ ($i=1,\dots ,q$).  
Theorem~\ref{opt.dual.set} and the definition we use 
for L$\sp {\natural}$-convex sets immediately implies the following.

\begin{corollary}   \label{Pistar} 
The set $\Pi$ of optimal dual integral vectors $\pi$ 
in the min-max formula \eqref{(minmax)} of 
Theorem~{\rm \ref{minmax}}
is an {\rm L}$\sp{\natural}$-convex set.  
The unique smallest element of $\Pi$ (that is, the unique smallest
dual optimum) is $\pi\sp{*}$.  
\finbox 
\end{corollary}

It will be worth mentioning that {\rm L}$\sp{\natural}$-convexity
of the set of optimal dual integral vectors
is a general phenomenon that is true in separable convex function minimization
on an M-convex set; see Section 5 of \cite{FM19partII}.
Indeed, this is a consequence of conjugacy between M-convexity and L-convexity.
It is also known that every {\rm L}$\sp{\natural}$-convex set has a 
unique smallest (and a unique largest) element.

\section{Continuous versus discrete}
\label{SCcontRdiscZ}

In this paper, we have concentrated on 
discrete decreasing minimization whose
continuous counterpart had been investigated earlier.  
In this section, we briefly look at the relationship 
between the continuous (fractional) and discrete (integral) settings.

First of all, there is a fundamental difference between the problems of
finding a dec-min element of a base-polyhedron and 
that of an M-convex set
(the set of integral elements of an integral base-polyhedron).  
In the  former case (investigated by Fujishige \cite{Fujishige80,Fujishigebook}), 
there alway exists a single, unique dec-min element, 
while in the latter, the dec-min elements
of an M-convex set have an elegant matroidal structure.  
Namely, Theorem~\ref{matroid-eltolt} shows that 
the set of dec-min elements of an M-convex set 
arises from the bases of a matroid by translating their
incidence vectors with an integral vector.

In spite of such fundamental difference,
a dec-min element can be characterized, in either case,
 as a square-sum minimizer.
In the continuous case, it is well-known \cite{Fujishige80,Fujishigebook} 
that the unique dec-min element of a base-polyhedron 
coincides with the minimum norm point of the base-polyhedron, 
whereas, in the discrete case,
Corollary \ref{COdecmin=squaresum} shows that the square-sum minimizers
are exactly the dec-min elements of an M-convex set. 
Furthermore,  a dec-min element can be characterized by 
a symmetric strictly convex function,
which is stated in Theorem~\ref{decmin=Phi=nonsep} 
for the discrete case
as a discrete counterpart of a result of Maruyama \cite{Mar78} 
for the continuous case.  
See also Nagano \cite[Corollary 13]{Nag07}.  
Symmetric convex function minimization is studied, mainly
for the continuous case, in the literature of majorization \cite{AS18,MOA11}.

In the following, we show links between 
the continuous and discrete versions of decreasing minimization
by considering an integral base-polyhedron $B$ and 
the asociated M-convex set $\odotZ{B}$.
The following theorems, given in \cite{FM19partII}, formalize 
the intuitive feeling that the minimum norm points (dec-min elements) 
of an M-convex set $\odotZ{B}$ 
and the unique minimum norm point (dec-min element) 
of the corresponding base-polyhedron $B$
are `close' to each other.

\begin{theorem}[{\cite[Theorem 6.6]{FM19partII}}] \label{THdecminproxS}
Let $m_{\RR}$ be the minimum norm point of an integral base-polyhedron $B$.
Then every dec-min element $m$ of $\odotZ{B}$ satisfies
$\left\lfloor m_{\RR} \right\rfloor \leq m \leq  \left\lceil m_{\RR} \right\rceil$.
\finbox
\end{theorem}

\begin{theorem}[{\cite[Theorem 6.7]{FM19partII}}] \label{THmnormconvcombdmS1}
The minimum norm point of an integral base-polyhedron $B$
can be represented as a convex combination of
the dec-min elements of $\odotZ{B}$.
\finbox
\end{theorem}

Since the (unique) minimum norm point of base-polyhedron $B$ is the
(unique) dec-min element of $B$, 
Theorem \ref{THmnormconvcombdmS1} can be reformulated in the following equivalent form. 

\begin{theorem} \label{THmnormconvcombdmS}
The (unique) dec-min element of an integral base-polyhedron $B$
can be represented as a convex combination of
the dec-min elements of $\odotZ{B}$.
\finbox
\end{theorem}

This result looks quite natural and even straightforward.
However, one has to be cautious with such a naive intuition.
In fact, the analogous statement 
fails to hold for an M$_{2}$-convex set 
(the intersection of two M-convex sets),
as is demonstrated in the following example.

\begin{example} \rm \label{EXconvhullRZ}
Consider the following two M-convex sets:  
\begin{align*} 
\odotZ{B_{1}} &= \{(1, 0, 0, 0), \ (0, -1, 1, 1), \ (1, -1, 1, 0), \ (0, 0, 0, 1)\}, 
\\ 
\odotZ{B_{2}} &= \{(1, 0, 0, 0), (0, -1, 1, 1), (1, \-1, 0, 1), (0, 0, 1, 0)\}.
\end{align*}
In their intersection 
$\odotZ{B_{1}} \cap \odotZ{B_{2}} = \{(1, 0, 0, 0), \ (0, -1, 1, 1)\}$,
the element $z = (1, 0, 0, 0)$ is the unique dec-min element.
In the continuous version,
$B_{1} \cap B_{2}$ is the line segment connecting $(1, 0, 0, 0)$ and $(0, -1, 1, 1)$.
The middle point
$x = (1/2,-1/2, 1/2, 1/2)$
is the unique dec-min element of $B_{1} \cap B_{2}$.
(Note that $x$ is decreasingly smaller than $z$.)
We cannot represent $x$ as a convex combination of a single element $z$.
\finbox
\end{example}

Not only the dec-min elements 
for $B$ and $\odotZ{B}$ 
are related as above,
but the `dual objects' (namely, the chains and partitions)
for $B$ and $\odotZ{B}$ are related as follows.
Actually, 
Theorems \ref{THdecminproxS} and \ref{THmnormconvcombdmS1} 
are proved in \cite{FM19partII}
 on the basis of this relationship between the `dual objects.' 
Recall that the principal partition is the continuous counterpart of 
the canonical partition, while the critical values are
the continuous counterpart of essential values
(see  \cite[Section 7.2]{Fujishigebook}, \cite{Fuj09bonn}
for notions related to the principal partition).

\begin{theorem}[{\cite[Theorem 6.5]{FM19partII}}]  \label{THrelRZpartition}
Let $B$ be an integral base-polyhedron.

\noindent
{\rm (1)}
An integer $\beta$ is an essential value of $\odotZ{B}$ 
if and only if 
there exists a critical value $\lambda$ of $B$ 
satisfying $\beta  \geq \lambda > \beta -1$.

\noindent
{\rm (2)}
The essential values $\beta_{1} > \beta_{2} > \cdots > \beta_{q}$
of $\odotZ{B}$ are obtained 
from the critical values 
$\lambda_{1} > \lambda_{2} >  \cdots >  \lambda_{r}$
of $B$
as the distinct members of the rounded-up integers
$\lceil \lambda_{1} \rceil \geq \lceil \lambda_{2} \rceil \geq \cdots \geq \lceil \lambda_{r} \rceil$.

\noindent
{\rm (3)}
The canonical partition 
$\{ S_{1}, S_{2}, \ldots, S_{q} \}$
of $\odotZ{B}$
is obtained from the principal partition 
$\{ \hat S_{1}, \hat S_{2}, \ldots, \hat S_{r} \}$
of $B$ as an aggregation as 
\[ 
 S_{j} = \bigcup_{i \in I(j)} \hat S_{i} 
\qquad (j=1,2,\ldots, q),
\] 
where $I(j) = \{ i :  \lceil \lambda_{i} \rceil  = \beta_{j} \}$ for $j=1,2,\ldots, q$.

\noindent
{\rm (4)}
The canonical chain
$\{ C_{j} \}$
of $\odotZ{B}$
is a subchain of the principal chain  $\{ \hat C_{i} \}$ of $B$,
which is given by $C_{j} = \hat C_{i}$ 
with the maximum index $i$ in $I(j)$.
\finbox
\end{theorem}

\section{Conclusion}
\label{SCconcl}

The present work is the first member of a series of papers concerning
discrete decreasing minimization.  
In the companion paper \cite{FM19partB} 
we give a strongly polynomial algorithm for finding a dec-min element of an M-convex set 
and discuss applications of discrete decreasing minimization 
to the `background problems' 
mentioned in Section~\ref{SCbackprob}.

While the present framework of 
decreasing minimization on an M-convex set
is effective for a fairly wide class of graph orientation problems \cite{FM19partB},
there are other important graph orientation problems that do not
fit in this framework.
For example, for strong orientations of mixed graphs, 
dec-min orientations and inc-max orientations do not coincide.  
The reason behind this phenomenon is that the set of in-degree vectors of 
strong orientations of a mixed graph is not an M-convex set anymore.  
It is, in fact, the intersection of two M-convex sets.
By investigating the decreasing minimization problem 
over the intersection of two M-convex sets
we can solve a broader class of graph orientation problems,
see \cite{FM20partD}.

Decreasing minimization on an M-convex set
contains the integer version of Megiddo's problem \cite{Megiddo74} 
of finding a maximum flow that is
\lq lexicographically optimal\rq \ on the set of edges 
leaving the source node.
In \cite{FM19partC} this problem is generalized  
to the problem of finding an integral feasible flow 
that is decreasing minimal on an arbitrarily specified subset of edges. 
The structure of decreasingly minimal integral feasible flows is clarified and 
a strongly polynomial algorithm for finding such a dec-min flow is developed.
A further generalization to integral submodular flows 
is reported in \cite{FM20partD}.

\paragraph{Acknowledgement} 
We are grateful to the six authors of
the paper by Borradaile et al.~\cite{BIMOWZ} because that work
triggered the present research (and this is so even if we realized
later that there had been several related works).  
We thank S. Fujishige and S. Iwata for discussion about the history of convex
minimization over base-polyhedra.  We also thank A. J\"uttner 
and T. Maehara for illuminating the essence of the Newton--Dinkelbach algorithm.  
J. Tapolcai kindly draw our attention to engineering
applications in resource allocation.  
Z. Kir\'aly played a similar role by finding
an article which pointed to a work of Levin and Onn on 
decreasingly minimal optimization in matroid theory.  
We are also grateful to M. Kov\'acs for drawing our attention to some
important papers in the literature concerning fair resource allocation problems.  
Special thanks are due to T. Migler for her continuous availability to answer our questions
concerning the paper \cite{BIMOWZ} and the work by Borradaile, Migler, and Wilfong \cite{BMW}, 
which paper was also a prime driving force in our investigations.
We are grateful to 
B. Shepherd and K. B{\'e}rczi
for their advice that led to restructuring our presentation appropriately. 
We are also grateful to an anonymous referee of the paper whose
strategic suggestions were particularly important to shape the final
form of our work.
This research was supported through the program ``Research in Pairs''
by the Mathematisches Forschungsinstitut Oberwolfach in 2019.
The two weeks we could spend at Oberwolfach provided 
an exceptional opportunity to conduct particularly intensive research.
The research was partially supported by the
National Research, Development and Innovation Fund of Hungary
(FK\_18) -- No. NKFI-128673,
and by CREST, JST, Grant Number JPMJCR14D2, Japan, 
and JSPS KAKENHI Grant Numbers JP26280004, JP20K11697. 






\end{document}